\documentclass[preprint,12pt]{elsarticle}

\usepackage{bbm}
\usepackage{xcolor}
\usepackage{hyperref}
\usepackage[a4paper, tmargin=1.5in, bmargin=1in, lmargin=1in, rmargin=1in, headheight=13.6pt]{geometry}
\usepackage{graphicx,epsfig}
\usepackage{enumerate}
\usepackage{amssymb}
\usepackage{amsthm}
\usepackage{amsmath}
\usepackage{centernot}
\usepackage{times}
\usepackage{mathrsfs}
\usepackage[utf8]{inputenc}
\usepackage{subfigure}
\usepackage{mathtools}
\usepackage{appendix}
\usepackage{wrapfig}
\usepackage{setspace}
\usepackage{booktabs}
\usepackage{lscape}
\usepackage{algorithm,algorithmic}
\usepackage{bm}
\usepackage{xfrac}
\usepackage{multirow}
\usepackage{longtable}
\usepackage{xcolor}
\usepackage{todonotes}
\usepackage{graphicx,epsfig}

\usepackage{wrapfig}
\usepackage{graphicx,epsfig}
\usepackage{xfrac}
\usepackage{caption}
\usepackage{times}
\usepackage{epstopdf}
\usepackage{amsmath}
\usepackage{paralist}
\usepackage{graphics}
\usepackage{epsfig}
\usepackage{graphicx}
\usepackage{epstopdf}
\usepackage{amsmath}

\usepackage{cleveref}
\hypersetup{
    colorlinks=true,
    linkcolor=blue,
    filecolor=magenta,
    urlcolor=cyan,
    }
\crefformat{section}{\S#2#1#3}
\crefformat{subsection}{\S#2#1#3}
\crefformat{subsubsection}{\S#2#1#3}
\crefrangeformat{section}{\S\S#3#1#4 to~#5#2#6}
\crefmultiformat{section}{\S\S#2#1#3}{ and~#2#1#3}{, #2#1#3}{ and~#2#1#3}

\newtheorem{theorem}{Theorem}[section]
\theoremstyle{definition}

\newtheorem{lemma}{Lemma}[section]

\newcommand{\sgn}{\operatorname{sgn}}

\newcommand{\decreasesto}[0]{\searrow}

\let\SavedIndent\indent
\protected\def\indent{%
  \begingroup
    \parindent=\the\parindent
    \SavedIndent
  \endgroup
}
\usepackage{natbib}

\makeatletter
\newcommand{\thickhline}{%
	\noalign {\ifnum 0=`}\fi \hrule height 1pt
	\futurelet \reserved@a \@xhline
}
\numberwithin{equation}{section}

\biboptions{sort&compress}
\begin{document}

\begin{frontmatter}
\title{\textbf{Mass conserving global solutions for the  nonlinear collision-induced fragmentation model with a singular kernel}}
%
%
\author[iitkgp_math]{Debdulal Ghosh\corref{cor1}}

\author[iitkgp_math]{Jayanta Paul}

\author[iitkgp_math]{Jitendra Kumar}
\address[iitkgp_math]{Department of Mathematics,
              Indian Institute of Technology Kharagpur,
              West Bengal 721302, India}
\cortext[cor1]{Corresponding author: debdulal.email@gmail.com (Debdulal Ghosh)}

\begin{abstract}This article is devoted to the study of existence of a mass conserving global solution for the collision-induced nonlinear fragmentation model which arises in particulate processes, with the following type of collision kernel:
\[C(x,y)~\le~k_1 \frac{(1 + x)^\nu (1 + y)^\nu}{\left(xy\right)^\sigma},\]
for all ~$x, y \in (0,\infty)$, where $k_1$ is a positive constant, $\sigma \in \left[0,\tfrac{1}{2}\right]$ and $\nu \in [0, 1]$. The above-mentioned form includes many practical oriented kernels of both \emph{singular} and \emph{non-singular} types. The singularity of the unbounded collision kernel at coordinate axes extends the previous existence result of Paul and Kumar [Mathematical Methods in the Applied Sciences 41 (7) (2018) 2715–2732 (\href{https://doi.org/10.1002/mma.4775}{doi:10.1002/mma.4775})] and also exhibits at most quadratic growth at infinity.
Finally, uniqueness of solution is also investigated for pure singular collision rate, i.e., for ~ $\nu=0$.
\end{abstract}
\begin{keyword} Nonlinear fragmentation phenomenon, Singular kernels, Global existence, Mass conservation, Uniqueness.\\
\vspace{0.5cm}
\end{keyword}
\end{frontmatter}

\section{Introduction}
\noindent \textcolor{black}{Reversible aggregation or Fragmentation (also, known as breakage) is the kinetic process of breaking up clusters by collision (nonlinear) or by cluster properties and external forces (linear). The significant kinetic process, fragmentation, has a diverse influences on many physical \cite{krapivsky2003shattering,brilliantov2015size,vyas2021collisional} as well as industrial processes \cite{yu2021particle,cohn2021dissociation,zhou2022experimental} and pharmaceutical research \cite{arzi2018electrohydrodynamic,nabizadeh2021life,cabiscol2021bi}. In this article, we focus on the nonlinear or collision-induced fragmentation as the theoretical and mathematical aspects of these models are not studied extensively like the linear one \cite{ziff1985kinetics,ramkrishna2000population,banasiak2014existence, breschi2017note}}.
\textcolor{black}{The mathematical modeling of the collision-induced fragmentation phenomenon has been developed by Cheng and Redner \cite{cheng1988scaling}, which has the following nonlinear integro-partial differential form}:
\begin{align}\label{main_eq}
\frac{\partial g(t, x)}{\partial t} =  \int_0^\infty \int_x^\infty C(y, z) F(x, y| z) g(t, y) g(t, z) dy dz
                                       - g(t, x) \int_0^\infty C(x, y) g(t, y) dy
\end{align}
\textcolor{black}{equipped with} initial data
\begin{align}\label{ini_cond}
g(0, x) = g_0(x), ~ x \in \mathbb{R}_+ = (0, \infty).
\end{align}
The temporal concentration of clusters \textcolor{black}{at time $t$ of mass or size or another physical characteristic (like density, volume, enthalpy, etc.)} between $x$ and $x+dx$ is denoted by $g(t, x)$. In equation (\ref{main_eq}), the \textcolor{black}{collisional} kernel $C(x, y)$ depicts the rate of effective \textcolor{black}{impact} for two clusters of \textcolor{black}{masses or} sizes $x$ and $y$, respectively. The \textcolor{black}{fragmentation} rate for forming an $x$-size cluster from a $y$-size cluster due to collision with a $z$-size particle is $F(x, y| z)$.
The collisional kernel $C(x, y)$ is symmetric \textcolor{black}{with reference to} size $x$ and $y$, i.e.,
\begin{align}
C(x, y) = C(y, x);~ \forall~ x>0, y>0 .
\end{align}
Concerning the \textcolor{black}{fragmentation rate} $F(x, y| z)$, it is a diversification of the \textcolor{black}{fragmentation kernel associated} in the linear \textcolor{black}{fragmentation model credible with an additional variable} $z$. The \textcolor{black}{fragmentation kernel $F(x, y| z)$ has the following characteristics}:
\begin{align}\label{prop_1_b(x,y;z)}
\int_0^y x F(x, y| z) dx = y; ~~ \text{ and } F(x, y| z) = 0; ~~\forall ~ x > y >0,
\end{align}
\begin{align}\label{prop_2_b(x,y;z)}
\int_0^y F(x, y| z) dx = \theta(y, z) < \infty.
\end{align}
The properties of (\ref{prop_1_b(x,y;z)}) indicate that the total \textcolor{black}{volume or mass of daughter} clusters \textcolor{black}{disintegrate} from the $y$ size \textcolor{black}{mother} cluster is $y$. In (\ref{prop_2_b(x,y;z)}), $\theta(y, z)$ represents the number of clusters due to the \textcolor{black}{split-up of a single $y$ size particle} on \textcolor{black}{successful collision for fragmentation} with a cluster of size $z$.
\textcolor{black}{In the preceding model equation (\ref{main_eq}):}
\begin{enumerate}[1.]
\item The first \textcolor{black}{integral} of \textcolor{black}{(\ref{main_eq})} is the formation of $x$ size cluster for the collision of $y$ and $z$ size clusters \textcolor{black}{with a precise} \textcolor{black}{fragmentation} rate $F(x, y| z)$.
\item The second \textcolor{black}{integral of (\ref{main_eq}) is the loss of $x$ size cluster owing to the collision of $x$ and $y$ size clusters.}
\end{enumerate}
\textcolor{black}{In this context, the mathematical moments of the system describe many physical characteristics, e.g., number, mass, volume, energy, etc.}
The $k$-th order moment is denoted by the following integral
\begin{align}
N_k(t) = \int_0^\infty x^k g(t, x) dx.
\end{align}
\textcolor{black}{In specific, the zeroth ($k=0$) and first ($k=1$) moments describe the aggregate
number and mass of the system. Also, the second moment($k=2$) illustrates the energy of the system.} A brief discussion and modeling of collisional fragmentation phenomenon can be found in \cite{safronov1972evolution,wilkins1982geometrical,cheng1988scaling,laurenccot2001discrete}
\section{\textcolor{black}{State of the art}}
The study of nonlinear collisional fragmentation model has already been the subject in the scientific community. However, the literature contains limited evidence on the study due to nonlinear behaviour and collision characteristics for the problem.
Cheng and Redner \cite{cheng1988scaling} have analyzed the scaling theory for the linear as well as nonlinear fragmentation process. They have shown scaled cluster-size distribution $\phi(x)$ for linear fragmentation, decrease to $x^{- 2} e^{- a x^\xi}$ as $x$ tends to infinity, ~$\xi$ is the homogeneity index. For small $x, \phi(x) \sim \exp( -a \ln^2 x)$ and a power-law form without the cut-off.
The scaling Ansatz for the cluster-size dynamics represent as $g(t,x) s^{-2} \phi(x/s),$ where $s$ is the conventional cluster characteristic (size, mass, volume, etc) . In comparing to the scaling form of the cluster-size dynamics for collisional fragmentation, they obtain
\[
 \phi(x)\sim
  \begin{cases}
    \exp(- x^{\xi/ 2})/ x^2\text{ both \textcolor{black}{the clusters break}, } \\
   \exp(- x^{\xi})/ x^2, \text{ larger \textcolor{black}{one breaks}, }\\
   x^{- (1 + \xi)}, \text{smaller \textcolor{black}{cluster breaks}, }
  \end{cases}
\]
for $x$ tends to infinity.
%
%
Cheng and Redner \cite{cheng1990kinetics} have also studied shattering transition, i.e., \textcolor{black}{volume} is lost due to the dust stage of sufficiently small size \textcolor{black}{clusters}. They also analyzed the dynamics of the collision-induced nonlinear fragmentation phenomenon. The asymptotic behavior of the time-dependent solution is investigated for the following collision results: (1) both clusters splintering into equal segments (2) only the larger cluster splintering in two (3) only the smaller cluster splintering. They have used scaling theory to obtain their result. Considering the linear fragmentation model, Banasiak and Lamb \cite{banasiak2014existence} have discussed the existence moments by semigroup operator theory. Furthermore, mathematical existence results of solution for discrete collisional fragmentation model with coagulation equations is studied by {Lauren{\c{c}}ot and Wrzosek} \cite{laurenccot2001discrete} in the following form
\begin{align} \label{laurencot1}
    \frac{dg_i}{dt}=&\frac{1}{2}\sum_{i=1}^{i-1}p_{j,i-j}k_{j,i-j}g_jg_{i-j}-\sum_{j=1}^{\infty}k_{i,j}g_ig_j \notag\\
    &+\frac{1}{2}\sum_{j=i+1}^{\infty}\sum_{k=j-1}^{\infty}(1-p_{j-k,k})f^i_{j-k,k}g_{j-k}g_k, ~i \geq 1
\end{align}
\begin{equation}\label{laurencot2}
    g_i(0)=g_i^0 ~i \geq 1
\end{equation}
In the preceding model \eqref{laurencot1}-\eqref{laurencot2},
\begin{itemize}
\item $k_{i,j}$ is the coagulation rate of two small size clusters of $i$ and $j$ to $i+j$ with probability $p_{i,j}$
\item \{$f^i_{j,k},~{i=1,..., j+k-1}\}$ is the distribution of fragments with probability of fragmentation  $1-p_{i,j}$.
\end{itemize}
In contrast to the continuous nonlinear model \eqref{main_eq}-\eqref{ini_cond}, here exist possible transfer of matter in fragmentation event. Mass conservation in fragmentation event indicates
\begin{equation}
    \sum_{i=1}^{j+k-1}if^i_{j,k}=j+k
\end{equation}
Now, the continuous Cheng and Redner \cite{cheng1988scaling} model \eqref{main_eq}-\eqref{ini_cond} is a generalised version of \eqref{laurencot1}-\eqref{laurencot2} for
$p_{i,j}=0, f^{i}_{j,k}=\chi_{[i,\infty)}(j)F^{i}_{j,k}+\chi_{[i,\infty)}(k)^{i}_{j,k}$, where $F^{i}_{j,k}$ is the corresponding discrete fragmentation rate of \eqref{main_eq}-\eqref{ini_cond}.
In this study, the authors have shown existence, density conservation and uniqueness of classical solution for unbounded non-singular kernels.
{Lauren{\c{c}}ot and Wrzosek} \cite{laurenccot2004time} have also studied the large time behavior of the corresponding Becker--D{\"o}ring-type of {coagulation model with collisional fragmentation process}. A specific case of coupled coagulation and singular kernel collisional breakage model is considered for the existence of a weak solution in \cite{barik2021existence,giri2021weak}. In another instance, the authors in \cite{cheng1990kinetics} have considered three models for collisional fragmentation of the form
\newline Model I: Both colliding clusters break into two by the process
\begin{align*}
& C(x,y)=x^{\xi/2}y^{\xi/2} \\
& F(x,y|z) =2 \delta(x-y/2)
\end{align*}
Model II: Just the bigger of the two impacting clusters break as
\begin{align*}
& C(x,y)=  \begin{cases}
   x^{\xi}  & \text{if } y\leq x \\
   y^{\xi} & \text{otherwise}.
  \end{cases} \\
& F(x,y|z) =\begin{cases}
   2 \delta(x-y/2)  & \text{if } z\leq y \\
   \delta(x-y) & \text{otherwise}.
  \end{cases}
\end{align*}
Model III: Only the smaller cluster breaks as
\begin{align*}
& C(x,y)=  \begin{cases}
   x^{\xi}  & \text{if } x\leq y \\
   y^{\xi} & \text{otherwise}.
  \end{cases} \\
& F(x,y|z) =\begin{cases}
   2 \delta(x-y/2)  & \text{if } y\leq z \\
   \delta(x-y) & \text{otherwise}.
  \end{cases}
\end{align*}
In the above models, the authors have shown all possible cases of scaling solutions for real $\xi$. Clearly, $\xi <0$ is nothing but the singular collisional kernel, which is also studied in \cite{cheng1990kinetics}. Analytical solutions, as well as self-similar analysis of the model problem, are discussed in \cite{kostoglou2000study} for non-singular collision rate. Although the theoretical aspects of collisional fragmentation have not been explored extensively, the corresponding discrete model of the problem has been studied by Lauren{\c{c}}ot and Wrzosek \cite{laurenccot2001discrete} for existence-uniqueness and mass conservation together with large time behavior for non-singular collisional kernels. Ernst and Pagonabarraga \cite{ernst2007nonlinear} have studied the scaling solution and shattering transition of the model for both singular and non-singular collision kernels.
Recently, Paul and Kumar \cite{paul2018existence} have derived the existence-uniqueness and mass conservation result by considering non-singular collision rate and singular fragmentation kernel of the form
\begin{equation*}
C(x,y) \leq \chi_0(1+x)^{\gamma}(1+y)^{\gamma}; ~F(x,y|z) \leq \chi_1\frac{1}{y^\beta}, ~\gamma,\beta> 0,
\end{equation*}
where $\chi_0$ and $\chi_1$ are positive constants. In their study \cite{paul2018existence}, the authors have shown the existence-uniqueness of mass conserving continuous solution in the space $G_a^+(T)$, for $a>0, T> 0$,~ where,
\begin{align*}
G_a^+(T):=\{\text{All non-negative continuous function}~ g :\sup_{0\leq t \leq T}\int_{0}^{\infty}x^a g(t,x)dx < \infty\}
\end{align*}
\par
Over the years many researchers have studied the nonlinear behaviour of colisional fragmentation model, e.g., scaling solutions \cite{cheng1988scaling, cheng1990kinetics,ernst2007nonlinear}, shattering behaviour \cite{krapivsky2003shattering,kostoglou2006study}, existence-uniqueness and well-posedness of weak solution \cite{giri2021existence}, Monte Carlo (Direct simulation) algorithm \cite{pagonabarraga2009collision}, discontinuous Galerkin scheme \cite{lombart2021grain}, etc.
However, to the best of authors knowledge the existence-uniqueness result of a global mass preserving continuous solution for singular kernel collision rate is not studied yet. The study of singular kernels are practically meaningful and  shed new light on known nonlinear phenomena \cite{saha2015singular,camejo2015singular,niethammer2016self,lamb2019continuous,ghosh2020existenceCAM,giri2021weak}. Physically, the singular kernel describes vast collision rate for smaller size clusters. In our study, the examined singular kernels  include a large class of physically meaningful kernels, e.g., Smoluchowski's Brownian diffusion kernel \cite{smoluchowski1917experiment}, Kapur kernel \cite{kapur1972kinetics}, Velocity (non-linear) profile kernel \cite{shiloh1973coalescence}, Equipartition (in granulation) kinetic energy kernel \cite{hounslow1998population}, Friedlander (in aerosol dynamics) \cite{friedlander2000smoke}, Ding et al. kernel \cite{ding2006population}, etc. This circumstance motivated the present study.

\section{Existence of solution}\label{sect_existence}
\textcolor{black}{In order to address the existence of a solution of (\ref{main_eq}), we construct the kernels $C_n$ in the following truncated form}
\[
 C_n(x, y)
  \begin{cases}
  \le C(x, y) & \text{if } (x, y) \in(0, \infty) \times (0, \infty) \setminus [\tfrac{1}{n}, n] \times [\tfrac{1}{n}, n] \\
  \decreasesto 0 & \text{in a finite range of }(0, \infty) \times (0, \infty) \setminus  [\tfrac{1}{n}, n] \times [\tfrac{1}{n}, n]\\
   = C(x, y)   & \text{if } (x, y)\in[\tfrac{1}{n}, n] \times [\tfrac{1}{n}, n].
  \end{cases}
\]
For the above "cutoff" kernels $C_n$ with solutions designated as $g^n$, the corresponding equation is
\begin{align}\label{truncated_main_eq}
\frac{\partial g^n(t, x)}{\partial t} =  \int_{\tfrac{1}{n}}^n \int_x^n C(y, z) F(x, y| z) g^n(t, y) g^n(t, z) dy dz
                                       - g^n(t, x) \int_{\tfrac{1}{n}}^n C(x, y) g^n(t, y) dy,
\end{align}
\textcolor{black}{Also, we truncate the initial data in the following approach}
\begin{align}\label{truncated_ini_cond}
 g^n(0, x) = g_0^n(x)=
  \begin{cases}
    g_0(x)   & \text{if } x\in[\tfrac{1}{n}, n] \\
    0 & \text{if } x \in(0, \infty) \setminus [\tfrac{1}{n}, n].
     \end{cases}
\end{align}
In the next, we define the strip $\mathcal{P}(T, X_1, X_2)$ in the following manner
\begin{align}
\mathcal{P}(T, X_1, X_2) = \{(t, x): t \in [0, T], ~ 0 < X_1 \le x \le X_2\},
\end{align}
where $T$ and $ X_1, X_2$ are finite numbers.
Let $T > 0$ be any given number and $0 < r \leq 1$. For a $\lambda > 0$, suppose $\Omega_{\lambda, r}(T)$ be the space of continuous functions $g$ in $\mathbb{R}_+ \times [0,T]$ with the norm
\[ \|g\|_{\lambda, r} \coloneqq \sup_{0 \le t \le T} \int_0^\infty \left( \exp(\lambda (1 + x)) + \tfrac{\exp(2 \lambda)}{x^r}\right)  |g(t,x)| dx. \]
With the help of $\Omega_{\lambda, r}(T)$, we define another function space $\Omega_{.,r}(T)$ by $$\Omega_{.,r}(T) \coloneqq \bigcup_{\lambda > 0} \Omega_{\lambda, r}(T).$$
Define $\Omega_{\lambda, r}^n(T)$ is the compactly supported continuous function $g$ over $[0, T] \times [\tfrac{1}{n}, n]$ with the norm
\[ \|g\|_{\lambda, r}^n \coloneqq \sup_{0 \le t \le T} \int_{\tfrac{1}{n}}^n\left( \exp(\lambda (1 + x)) + \tfrac{\exp(2 \lambda)}{x^r}\right)  |g(t,x)| dx < \infty. \]
Cones of non negative functions in $\Omega_{.,r}(T)$ is denoted by $\Omega_{.,r}^+(T)$.
Throughout the section, we take the following assumptions on the kernels.\\
\\
%
\textbf{\emph{Hypotheses:}}
\begin{enumerate}[(H1)]
\item The non-negative collision kernel $C(x, y)$ is \textcolor{black}{continuous} in $(0, \infty) \times (0, \infty)$.
\item The \textcolor{black}{non-negative fragmentation kernel} $F(x, y| z)$ is \textcolor{black}{continuous} in $(0, \infty) \times (0, \infty) \times (0, \infty)$.
\item The collision kernel satisfies ~$C(x,y) \le k_1~ \frac{(1 + x)^\nu (1 + y)^\nu}{\left(xy\right)^\sigma}$, $\forall ~x, y \in (0,\infty),$   where~~ $k_1 > 0$ is a constant, $\sigma \in \left[0,\tfrac{1}{2}\right]$ and $\nu \in [0, 1]$.
\item For all $0 < x < y,$ there exist real number $0 < \beta \le \sigma$, so that $F(x, y| z) \le \frac{k_2}{y^\beta}$, where $k_2$ is a positive constant.
\end{enumerate}
%
\begin{theorem}\textbf{(Local existence-uniqueness)}\label{thm_truncated}\\
Assume the collision rate $C(x, y)$ and the fragmentation rate $F(x, y| z)$ be non-negative and continuous in $(0, \infty) \times (0, \infty)$. In addition, let $C(x, y)$ is symmetric with respect to its arguments $x$ and $y$ in $(0, \infty) \times (0, \infty)$
and the initial data function satisfy $g_0 \in \Omega_{.,r_2}^+ (0)$. Then the truncated model (\ref{truncated_main_eq}) has an unique solution $g^n \in \Omega_{\lambda,r_2}^n (T)$ for each $n = 1, ~2, ~3,...$ and for $0 \leq t \leq T$, $x \in [\tfrac{1}{n}, n]$. \textcolor{black}{Moreover}, the conservation (mass) \textcolor{black}{property} holds, i.e.,
\begin{align}\label{truncated_mass_consv}
\int_{\tfrac{1}{n}}^n x g^n(t, x) dx = \int_{\tfrac{1}{n}}^n x g^n(0, x) dx, \text{ for } 0 \le t \le T.
\end{align}
\end{theorem}
\begin{proof}
The equation (\ref{truncated_main_eq}) can be reformulated in the following fashion:
\begin{align}\label{reform_truncated_main_eq}
\frac{\partial}{\partial t} [\exp(D(t, x, g^n)) g^n(t, x)] =  & \exp(D(t, x, g^n)) \notag\\
                                                              & \times \left[ \int_{\tfrac{1}{n}}^n \int_x^n C(y, z) F(x, y| z) g^n(t, y) g^n(t, z) dy dz \right]
\end{align}
supported by the initial datum
\begin{align}\label{truncated_ini_cond_1}
 g^n(0, x) = g_0^n(x)=
  \begin{cases}
   = g_0(x)   & \text{if } x\in[\tfrac{1}{n}, n] \\
   = 0 & \text{if } x \in(0, \infty) \setminus [\tfrac{1}{n}, n],
     \end{cases}
\end{align}
where
\begin{align}
D(t, x, g^n) = \int_0^t \int_{\tfrac{1}{n}}^n C(x, y) g^n(s, y) dy ds.
\end{align}
By integrating (\ref{reform_truncated_main_eq}) in the range $[0,~ t]$, we get
\begin{align}\label{iteration_f_n}
g^n(t, x) = \mathcal{C}(g^n)(t, x),
\end{align}
where
\begin{align}
\mathcal{C}(g^n)(t, x) = g_0^n & \exp(- D(t, x, g^n)) + \int_0^t  \exp\{-(D(t, x, g^n) - D(s, x, g^n))\}\notag\\
                     & \times \int_{\tfrac{1}{n}}^n \int_x^n C(y, z) F(x, y| z) g^n(t, y) g^n(t, z) dy dz ds.
\end{align}
\textcolor{black}{Now, we shall go through some important lemmas, in order to prove the theorem}. With the help of these lemmas, we show $\mathcal{C}$ has a fixed point by contraction mapping theorem in $[0, t_0]$.
Let us choose
\begin{align}
L = \|g_0\|_{\lambda, r} + \left(\frac{\exp(\lambda(1 + n))}{\lambda}+ \frac{\exp(2 \lambda) n^{1 - r}}{1 - r}\right) M^2 T \|g_0\|_1^2.
\end{align}
and \textcolor{black}{fix $t', t^{''} > 0$, so that}
\begin{align}\label{less_2}
\exp(2 t M L) \left(1 + 4 L K_1 t\left(\frac{\exp(\lambda(1 + n))}{\lambda}+ \frac{\exp(2 \lambda) n^{1 - r}}{1 - r}\right) \right) \le 2, \text{ for } 0 \le t \le t',
\end{align}
\begin{align}\label{k<1}
& \exp(t B M) \left( M t \|g_0\|_{\lambda, r}^n + k_1 k_2\left(\frac{\exp(\lambda(1 + n))}{\lambda}+ \frac{\exp(2 \lambda) n^{1 - r}}{1 - r}\right) (M t^2 B^2 + 2 B t)\right)  \notag\\
& < 1,~~ \text{ for } 0 \le t \le t^{''}.
\end{align}
Also, \textcolor{black}{choose} $t_0$ and $M$ as,
\begin{align}
\max \{ &\sup \{ F(x, y| z); x, z, y\in [\tfrac{1}{n}, n]\};
            \sup\{C(x,y); x, y\in [\tfrac{1}{n}, n] \}\}=M,\notag\\
            &~~~~~~\min\{t^{''}, t', T \}=t_0.
\end{align}
\begin{lemma}\label{aux_lemma}
For $g_1, g_2 \in \Omega_{\lambda, r}^n(t_0)$ and $t_0 \geq t \geq s \geq 0; ~ \tfrac{1}{n} \le x \le n$, \textcolor{black}{following result holds}
\begin{align}\label{ineq_F}
|F(s, t, x)| \le M (t - s) \|g_1 - g_2\|_{\lambda, r}^n \exp((t - s) B M),
\end{align}
where, $\exp\{- (D(t, x, g_1) -D(s, x, g_1))\} - \exp\{- (D(t, x, g_2) -D(s, x, g_2))\}=F(s, t, x)$ and $B = \max \{ \|g_1\|_{\lambda, r}^n, \|g_2\|_{\lambda, r}^n\}$.
\end{lemma}
\begin{proof}
\textcolor{black}{We consider}, $D(t, x, g_1) - D(s, x, g_1) \ge D(t, x,g_2) - D(s, x, g_2)$. Therefore
\begin{align}
|F(s, t, x)| = & - F(s, t, x)= \exp\{- (D(t, x, g_2) -D(s, x, g_2))\} \notag\\
               & \{ 1 - \exp\{- [D(t, x, g_1) -D(s, x, g_1) - (D(t, x, g_2) -D(s, x, g_2))]\} \}.
\end{align}
As $1 - \exp(- x) \le x$, for $x \ge 0$, we get
\begin{align}
|F(s, t, x)| & \le \exp\{- (D(t, x, g_2) -D(s, x, g_2))\} \notag\\
             &  \times \{ (D(t, x, g_1) -D(s, x, g_1)) - (D(t, x, g_2) -D(s, x, g_2))\}\notag\\
             &\le \exp\{- (D(t, x, g_2) -D(s, x, g_2))\} \int_s^t \int_{\tfrac{1}{n}}^n C(x, y) [g_1(\tau, y) -g_2(\tau, y)] dy d\tau \notag\\
             & \le \exp \left( \int_s^t \int_{\tfrac{1}{n}}^n C(x, y) g_2(\tau, y) dy d\tau\right)\notag\\
             & ~~~ ~~~\times M \int_s^t \int_{\tfrac{1}{n}}^n \left( \exp(\lambda(1 + y)) + \frac{\exp(2 \lambda)}{y^r} \right) |g_1 - g_2|(y, \tau) dy d\tau \notag\\
             & \le \exp(M B (t - s)) M \|g_1 - g_2\|_{\lambda, r}^n (t - s).
\end{align}
If $D(t, x, g_1) - D(s, x, g_1) \le D(t, x, g_2) - D(s, x, g_2)$, then the inequality (\ref{ineq_F}) can be proved by similar way.
\end{proof}
\begin{lemma}\label{C_into_itself}
The nonlinear integral operator $\mathcal{C}$ is a mapping from $\Xi$ to $\Xi$, where
\[\Xi = \{ g \in \mathcal{C}([0, t_0] \times [\tfrac{1}{n}, n]): \|g\|_{\lambda, r}^n \le 2L \}.\]
\end{lemma}
\begin{proof}
Let us assign $ t \in[0, t_0]$ and $\|g\|_{\lambda, r}^n \le 2L$, we get
\begin{align}
& \int_{\tfrac{1}{n}}^n  \left( \exp(\lambda(1 + x)) + \frac{\exp(2 \lambda)}{x^r}\right) |\mathcal{C}(g)(t, x)| dx \notag\\
                      & \le \underbrace{\int_{\tfrac{1}{n}}^n \left( \exp(\lambda(1 + x)) + \frac{\exp(2 \lambda)}{x^r}\right) g_0(x)\exp(- D(t, x, g)) dx}_{A_1}\notag\\
                      & ~~~  \underbrace{\begin{aligned} & + \int_{\tfrac{1}{n}}^n \int_0^t \left( \exp(\lambda(1 + x)) + \frac{\exp(2 \lambda)}{x^r}\right) \exp(E(t, x, g)) \notag\\
                                                         & ~~~    ~~~ \times \int_{\tfrac{1}{n}}^n \int_x^n C(y, z) F(x, y| z) |g(s, y)| |g(s, z)| dy dz ds dx,
                                                         \end{aligned}}_{A_2}
\end{align}
%
%
where, $[D(s, x, g) - D(t, x, g)]=E(t, x, g)$ and
\begin{align}
A_1 =   & \int_{\tfrac{1}{n}}^n \left( \exp(\lambda(1 + x)) + \frac{\exp(2 \lambda)}{x^r}\right) \exp(- D(t, x, g)) dx\notag\\
    \le & \int_{\tfrac{1}{n}}^n \left( \exp(\lambda(1 + x)) + \frac{\exp(2 \lambda)}{x^r}\right) \exp( D(t, x, g)) dx \notag\\
    \le & \exp(t M \|g\|_{\lambda, r}^n) \|g_0\|_{\lambda, r}^n.
\end{align}
%
Using Fubini's theorem and a change in the order of integration in $A_2$, we obtain
\begin{align}\label{A_2}
A_2  &  \le \int_{\tfrac{1}{n}}^n \int_0^t \left( \exp(\lambda(1 + x)) + \frac{\exp(2 \lambda)}{x^r}\right) \exp(D(t, x, g)) \notag\\
     &  \times \int_{\tfrac{1}{n}}^n \int_x^n C(y, z) F(x, y| z) |g(s, y)| |g(s, z)| dy dz ds dx, \notag\\
     & \le \exp(t M \|g\|_{\lambda, r}^n) \int_0^t \int_{\tfrac{1}{n}}^n   \int_{\tfrac{1}{n}}^n \int_x^n \left( \exp(\lambda(1 + x)) + \frac{\exp(2 \lambda)}{x^r}\right) \notag\\
     & ~~~ C(y, z) F(x, y| z) |g(s, y)| |g(s, z)| dy dz dx ds.
\end{align}
Now,
\begin{align}
& \int_{z = \tfrac{1}{n}}^n   \int_{x =\tfrac{1}{n}}^n \int_{y = x}^n \left( \exp(\lambda(1 + x)) + \frac{\exp(2 \lambda)}{x^r}\right)
C(y, z) F(x, y| z) |g(s, y)| |g(s, z)| dy  dx dz \notag\\
& = \int_{z = \tfrac{1}{n}}^n   \int_{y =\tfrac{1}{n}}^n \int_{x = \tfrac{1}{n}}^y \left( \exp(\lambda(1 + x)) + \frac{\exp(2 \lambda)}{x^r}\right)
C(y, z) k_2 y^{- \beta} |g(s, y)| |g(s, z)| dy dz dx \notag\\
& \le \left( \exp(\lambda(1 + n)) \lambda^{-1} + \frac{\exp(2 \lambda)n^{1 - r}}{1 - r}\right)  \notag\\
&   ~~~ ~~~ \int_{z = \tfrac{1}{n}}^n   \int_{y =\tfrac{1}{n}}^n k_1 k_2 \frac{(1 + y)^\nu (1 + z)^\nu}{(yz)^\sigma} y^{- \beta} |g(s, y)| |g(s, z)| dy dz dx \notag\\
& = k_1 k_2\left( \exp(\lambda(1 + n)) \lambda^{-1} + \frac{\exp(2 \lambda)n^{1 - r}}{1 - r}\right)  \notag\\
&   ~~~ ~~~  \int_{y =\tfrac{1}{n}}^n  \frac{(1 + y)^\nu}{y^{\sigma + \beta}} |g(s, y)| dy \int_{z = \tfrac{1}{n}}^n   \frac{(1 + z)^\nu}{(z)^\sigma}  |g(s, z)| dz.
\end{align}
Next we obtain
\begin{align}\label{singular_ineq_norm}
& \int_{y =\tfrac{1}{n}}^n  \frac{(1 + y)^\nu}{y^{\sigma + \beta}} |g(s, y)| dy \notag\\
& = \int_{y =\tfrac{1}{n}}^1  \frac{(1 + y)^\nu}{y^{\sigma + \beta}} |g(s, y)| dy + \int_{y =1}^n  \frac{(1 + y)^\nu}{y^{\sigma + \beta}} |g(s, y)| dy \notag\\
& = \int_{y =\tfrac{1}{n}}^1  \frac{\exp(\lambda (1 + y))}{y^{\sigma + \beta}} |g(s, y)| dy + \int_{y =1}^n  \frac{\exp(\lambda (1 + y))}{y^{\sigma + \beta}} |g(s, y)| dy \notag\\
& \le \int_{y =\tfrac{1}{n}}^1  \frac{\exp(2 \lambda)}{y^r} |g(s, y)| dy + \int_{y =1}^n  \exp(\lambda (1 + y)) |g(s, y)| dy \text{ as } \sigma + \beta \le r\notag\\
& = \int_{y =\tfrac{1}{n}}^n \left( \exp(\lambda (1 + y)) + \frac{\exp(2 \lambda)}{y^r} \right) |g(s, y)| dy\notag\\
& = \|g\|_{\lambda, r}^n.
\end{align}
Similarly,
\begin{align}
 \int_{z = \tfrac{1}{n}}^n   \frac{(1 + z)^\nu}{(z)^\sigma}  |g(s, z)| dz  \le \|g\|_{\lambda, r}^n.
\end{align}
Therefore, form (\ref{A_2}), we obtain
\begin{align}
A_2 \le \exp(t M \|g\|_{\lambda, r}^n)  k_1 k_2\left( \exp(\lambda(1 + n)) \lambda^{-1} + \frac{\exp(2 \lambda)n^{1 - r}}{1 - r}\right) (\|g\|_{\lambda, r}^n)^2 t.
\end{align}
So, using (\ref{less_2}), we get
\begin{align}
&\int_{\tfrac{1}{n}}^n  \left( \exp(\lambda(1 + x)) + \frac{\exp(2 \lambda)}{x^r}\right) |\mathcal{C}(g)(t, x)| dx \notag\\
                      & \le \exp(t M \|g\|_{\lambda, r}^n)  \left\{ \|g_0\|_{\lambda, r}^n +  k_1 k_2 t \left( \exp(\lambda(1 + n)) \lambda^{-1} + \frac{\exp(2 \lambda)n^{1 - r}}{1 - r}\right) (\|g\|_{\lambda, r}^n)^2 \right\} \notag\\
                      & \le \exp(t M 2 L)  \left\{ L +  k_1 k_2 t \left( \exp(\lambda(1 + n)) \lambda^{-1} + \frac{\exp(2 \lambda)n^{1 - r}}{1 - r}\right) 4L^2 \right\} \notag\\
                      & \le L ~~\underbrace{\exp(t M 2 L)  \left\{ 1 +  k_1 k_2 t \left( \exp(\lambda(1 + n)) \lambda^{-1} + \frac{\exp(2 \lambda)n^{1 - r}}{1 - r}\right) 4L \right\}}_{\le 2} \notag\\
                      & \le 2L.
\end{align}
\end{proof}
\begin{lemma}\label{C_contraction}
The integral operator $\mathcal{C}$ preserves contraction mapping property on $\Xi$, where
\[\Xi = \{ g \in \mathcal{C}([0, t_0] \times [\tfrac{1}{n}, n]): \|g\|_{\lambda, r}^n \le 2L \}.\]
\end{lemma}
\begin{proof}
Consider $f, g \in \Xi$ and recall $F$ from Lemma \ref{aux_lemma}. Therefore
\begin{align}
& \mathcal{C}(f)(t, x) - \mathcal{C}(g)(t, x) = g_0^n(x) F(0, t, x) \notag\\
& + \int_0^t F(t, s, x) \int_{\tfrac{1}{n}}^n \int_x^n C(y, z) F(x, y| z) g^n(s, y) g^n(s, z) dy dz ds \notag\\
& + \int_0^t \exp([D(s, x, g) - D(t, x, g)]) [\bar{B}(s, x, f) - \bar{B}(s, x, g)] ds \notag\\
\text{ and } &  \bar{B}(s, x, f) = \int_{\tfrac{1}{n}}^n \int_x^n C(y, z) F(x, y| z) g^n(s, y) g^n(s, z) dy dz.
\end{align}
Hence, we obtain
\begin{align}\label{B_123}
\int_{\tfrac{1}{n}}^n  \left( \exp(\lambda(1 + x)) + \frac{\exp(2 \lambda)}{x^r}\right) |\mathcal{C}(f)(t, x) - \mathcal{C}(g)(t, x)| dx \le \sum_{i=1}^{3} B_i,
\end{align}
where
\begin{align}
B_1 & = \int_{\tfrac{1}{n}}^n  \left( \exp(\lambda(1 + x)) + \frac{\exp(2 \lambda)}{x^r}\right) g_0^n(x) F(x, 0, t) dx \notag\\
    & \le M t \|f - g\|_{\lambda, r}^n \exp(t B M) \|g_0\|_{\lambda, r}^n,
\end{align}
the second estimate of (\ref{B_123}) becomes
\begin{align}
B_2 & = \int_0^t \int_{\tfrac{1}{n}}^n \int_{\tfrac{1}{n}}^n \int_x^n  \left( \exp(\lambda(1 + x)) + \frac{\exp(2 \lambda)}{x^r}\right) \notag\\
    & ~~~~~ ~~~~~ ~~~~~ ~~~~~ ~~~~~ F(t, s, x) C(y, z) F(x, y| z) f(s, y) f(s, z) dy dz dx ds \notag\\
    & \le M (t - s) \|f - g\|_{\lambda, r}^n \exp((t - s) B M) \int_0^t \int_{\tfrac{1}{n}}^n \int_{y = \tfrac{1}{n}}^n \int_{x = \tfrac{1}{n}}^y \notag\\
    &  ~~~~~ ~~~~~ \left( \exp(\lambda(1 + x)) + \frac{\exp(2 \lambda)}{x^r}\right)  C(y, z) F(x, y| z) f(s, y) f(s, z) dx dy dz  ds \notag\\
    & \le M (t - s) \|f - g\|_{\lambda, r}^n \exp((t - s) B M) \left( \exp(\lambda(1 + n)) \lambda^{-1} + \frac{\exp(2 \lambda)n^{1 - r}}{1 - r}\right) \notag\\
    & ~~~~~ ~~~~~ \int_0^t \int_{y = 0}^n \int_{z = \tfrac{1}{n}}^n   k_1 \frac{(1 + y)^\nu (1 + z)^\nu}{(yz)^\sigma} \frac{k_2}{y^\beta}  f(s, y) f(s, z) dy dz dx ds \notag\\
    & \le k_1 k_2 M (t - s) \|f - g\|_{\lambda, r}^n \exp((t - s) B M)  \notag\\
    &~~~~~ ~~~~~ \left( \exp(\lambda(1 + n)) \lambda^{-1} + \frac{\exp(2 \lambda)n^{1 - r}}{1 - r}\right) (\|f\|_{\lambda, r}^n)^2 t, \text{ using similar calculation in (\ref{singular_ineq_norm})}, \notag\\
    & \le k_1 k_2 M (t - s) \|f - g\|_{\lambda, r}^n \exp((t - s) B M)   \left( \exp(\lambda(1 + n)) \lambda^{-1} + \frac{\exp(2 \lambda)n^{1 - r}}{1 - r}\right) B^2 t,
\end{align}
and
\begin{align*}
B_3 & = \int_0^t \int_{\tfrac{1}{n}}^n  \left( \exp(\lambda(1 + x)) + \frac{\exp(2 \lambda)}{x^r}\right) \exp(-[D(t, x, g) - D(s, x, g)]) \notag\\
    & ~~~~~ ~~~~~ ~~~~~ [\bar{B}(s, x, f) - \bar{B}(s, x, g)] dx ds \notag\\
    & \le \exp(t M B) \int_0^t \int_{\tfrac{1}{n}}^n  \left( \exp(\lambda(1 + x)) + \frac{\exp(2 \lambda)}{x^r}\right) [\bar{B}(s, x, f) - \bar{B}(s, x, g)] dx ds \notag\\
    & \le \exp(t M B) \int_0^t \int_{\tfrac{1}{n}}^n \int_{\tfrac{1}{n}}^n \int_x^n  \left( \exp(\lambda(1 + x)) + \frac{\exp(2 \lambda)}{x^r}\right)\notag\\
    & ~~~~~ ~~~~~ ~~~~~ C(y, z) F(x, y| z) |f(s,y) f(s, z) - g(s, y) g(s, z)| dy dz dx ds \notag\\
        \end{align*}
    \begin{align}
    & \le \exp(t M B) \int_0^t \int_{z = \tfrac{1}{n}}^n \int_{y = \tfrac{1}{n}}^n \int_{x=\tfrac{1}{n}}^y  \left( \exp(\lambda(1 + x)) + \frac{\exp(2 \lambda)}{x^r}\right)\notag\\
    & ~~~~~ ~~~~~ ~~~~~ k_1 \frac{(1 + y)^\nu (1 + z)^\nu}{(yz)^\sigma} \frac{k_2}{y^\beta} |f(s,y) f(s, z) - g(s, y) g(s, z)| dx dy dz ds \notag\\
    & \le \exp(t M B) \left( \exp(\lambda(1 + n)) \lambda^{-1} + \frac{\exp(2 \lambda)n^{1 - r}}{1 - r}\right) \int_0^t \int_{z = \tfrac{1}{n}}^n \int_{y = \tfrac{1}{n}}^n   \notag\\
    & ~~~~~ ~~~~~ ~~~~~ k_1 \frac{(1 + y)^\nu (1 + z)^\nu}{(yz)^\sigma} \frac{k_2}{y^\beta} |f(s,y) f(s, z) - g(s, y) g(s, z)| dy dz ds \notag\\
    & \le \exp(t M B) \left( \exp(\lambda(1 + n)) \lambda^{-1} + \frac{\exp(2 \lambda)n^{1 - r}}{1 - r}\right) \int_0^t \int_{z = \tfrac{1}{n}}^n \int_{y = \tfrac{1}{n}}^n k_1 \frac{(1 + y)^\nu (1 + z)^\nu}{(yz)^\sigma}   \notag\\
    & ~~~~~ ~~~~~ \frac{k_2}{y^\beta} (g(s,y) |g(s, z) - f(s, z)| + f(s, z)|g(s, y)- f(s, y)|) dy dz ds. %
\end{align}
Now,
\begin{align}
& \int_{z = \tfrac{1}{n}}^n \int_{y = \tfrac{1}{n}}^n k_1 \frac{(1 + y)^\nu (1 + z)^\nu}{(yz)^\sigma} \frac{k_2}{y^\beta} g(s,y) |g(s, z) - f(s, z)|  dy dz \notag\\
& = k_1 k_2 \int_{y = \tfrac{1}{n}}^n \frac{\exp(\lambda (1 + y))}{y^{\sigma + \beta}} g(s, y) dy \int_{z = \tfrac{1}{n}}^n \frac{\exp(\lambda(1 + z))}{z^\sigma} |g(s, z) - f(s, z)| dz \notag\\
&\le k_1 k_2 \|g\|_{\lambda, r}^n \|g - f\|_{\lambda, r}^n \text{ using similar calculation in (\ref{singular_ineq_norm})},\notag\\
& \le k_1 k_2 B \|g - f\|_{\lambda, r}^n.
\end{align}
Similarly,
\begin{align}
& \int_{z = \tfrac{1}{n}}^n \int_{y = \tfrac{1}{n}}^n k_1 \frac{(1 + y)^\nu (1 + z)^\nu}{(yz)^\sigma} \frac{k_2}{y^\beta} f(s, z)|g(s, y)- f(s, y)|  dy dz \notag\\
& = k_1 k_2 \int_{y = \tfrac{1}{n}}^n \frac{\exp(\lambda (1 + y))}{y^{\sigma + \beta}} |g(s, y)- f(s, y)|  dy \int_{z = \tfrac{1}{n}}^n \frac{\exp(\lambda(1 + z))}{z^\sigma} f(s, z) dz \notag\\
&\le k_1 k_2  \|g - f\|_{\lambda, r}^n \|f\|_{\lambda, r}^n,  \text{ using similar calculation in (\ref{singular_ineq_norm})},\notag\\
& \le k_1 k_2 B \|g - f\|_{\lambda, r}^n.
\end{align}
Therefore,
\begin{align}
B_3 \le \exp(t M B) \left( \exp(\lambda(1 + n)) \lambda^{-1} + \frac{\exp(2 \lambda)n^{1 - r}}{1 - r}\right) 2 k_1 k_2 B t \|g - f\|_{\lambda, r}^n.
\end{align}
Hence from (\ref{B_123}), for $0 \le t \le t_0$, we obtain
\begin{align}
\|\mathcal{C}(f) - \mathcal{C}(g)\|_{\lambda, r}^n & = \int_{\tfrac{1}{n}}^n \left( \exp(\lambda(1 + x)) + \frac{\exp(2 \lambda)}{x^r}\right) |          \mathcal{C}(f) - \mathcal{C}(g)|(t, x) dx \notag\\
                                                   & \le \|f - g\|_{\lambda, r}^n \exp(t B M) \bigg[ M t \|g_0\|_{\lambda, r}^n  \notag\\
                                                   & ~~~~~ ~~~~~ + k_1 k_2 \left( \exp(\lambda(1 + n)) \lambda^{-1} + \frac{\exp(2 \lambda)n^{1 - r}}{1 - r}\right) (M B^2 t^2 + 2 B t) \bigg] \notag\\
                                                   & = k \|f - g\|_{\lambda, r}^n,
\end{align}
where
\begin{align}\label{contraction_k}
k  = \exp(t B M) \bigg[& M t \|g_0\|_{\lambda, r}^n + k_1 k_2 \left( \exp(\lambda(1 + n)) \lambda^{-1} + \frac{\exp(2 \lambda)n^{1 - r}}{1 - r}\right) \notag\\
                                                   & \times (M B^2 t^2 + 2 B t) \bigg].
\end{align}
By (\ref{k<1}), $k < 1$, therefore, we get the desired result.
\end{proof}
As a result of Lemma \ref{C_into_itself}, Lemma \ref{C_contraction}, and the Banach Fixed point theorem, we obtain a unique solution $g^n$ for $ 0 \le t \le t_0$. Furthermore, if we establish $g_0 = g_0^n$, $g_{\alpha} = \mathcal{C}(g_{\alpha - 1})$ for $\alpha = 1, 2, 3,...$, then $g_{\alpha} \rightarrow g^n$ in $\Omega_{\lambda, r}^+(t_0)$ as $\alpha \rightarrow \infty$, whereby (\ref{truncated_ini_cond_1}) and (\ref{iteration_f_n}), $g^n$ is \textcolor{black}{constructed} through positivity preservation. In the next, for $0 \le t \le t_0$ from (\ref{truncated_main_eq}), we can deduce that
\begin{align}\label{mass_consv_truncated_aux_eq}
\frac{d}{dt} \left(\int_{\tfrac{1}{n}}^n x g^n(t, x) dx \right) = & \int_{\tfrac{1}{n}}^n \int_{\tfrac{1}{n}}^n \int_x^n x C(y, z) F(x, y| z) g^n(t, y) g^n(t, z) dy dz dx\notag\\
                                                      &- \int_{\tfrac{1}{n}}^n g^n(t, x) \int_{\tfrac{1}{n}}^n x C(x, y) g^n(t, y) dy dx.
\end{align}
Due to the continuity of the kernels $C(x, y)$ and $F(x, y| z)$ on compact domain, the integrals of (\ref{mass_consv_truncated_aux_eq}) are finite. The first integration on r.h.s. of  (\ref{mass_consv_truncated_aux_eq}) becomes
%
%
\begin{align}
\int_{\tfrac{1}{n}}^n \int_{\tfrac{1}{n}}^n \int_x^n x C(y, z) F(x, y| z) g^n(t, y) g^n(t, z) dy dz dx
= \int_{\tfrac{1}{n}}^n \int_{\tfrac{1}{n}}^n   y C(y, z) g^n(t, y) g^n(t, z) dy dz.
\end{align}
Hence, we get
\[\frac{d}{dt} \left(\int_{\tfrac{1}{n}}^n x g^n(t, x) dx \right) = 0. \]
Therefore the result (\ref{truncated_mass_consv}) holds good. To enlarge the range of $t$ from $[0, t_0]$ to $[0, T]$, observe the following form
\[g^n(t, x) = \mathcal{C}_1(g^n)(t, x) - \mathcal{C}_1(g^n)(t_0, x),\]
whereby, the form of $\mathcal{C}_1$ is given below
\begin{align}
&\mathcal{C}_1(g^n)  = g^n(t_0) \exp(-D_1(t, x, g^n)) \notag\\
& + \int_{t_0}^t \exp(-(D_1(t, x, g^n) - D_1(s, x, g^n))) \int_{\tfrac{1}{n}}^n \int_x^n x C(y, z) F(x, y| z) g^n(t, y) g^n(t, z) dy dz ds,\\
&\text{ and } D_1(t, x, g) = \int_{t_0}^t \int_{\tfrac{1}{n}}^n C(x, y) g(s, y) dy ds.
\end{align}
Continuing as before, we can prove that there exists a non-negative unique solution on $[t_0, t_1]$, where $t_0 < t_1$. Time intervals can be extended by repeating the technique to $[0, T]$.
By integrating (\ref{truncated_main_eq}) multiplied by $\left( \exp(\lambda(1 + x)) + \frac{\exp(2 \lambda)}{x^r}\right)$, successionally w.r.t~ $t$ ~and~ $x$, we get
\begin{align}\label{to_show_f_n_in_Omega}
& \int_{\tfrac{1}{n}}^n \left( \exp(\lambda(1 + x)) + \frac{\exp(2 \lambda)}{x^r}\right) g^n(t, x) dx \notag\\
& =  \int_{\tfrac{1}{n}}^n \left( \exp(\lambda(1 + x)) + \frac{\exp(2 \lambda)}{x^r}\right) g_0(t, x) dx \notag\\
& + \int_{\tfrac{1}{n}}^n \int_0^t \int_{\tfrac{1}{n}}^n \int_x^n x  \left( \exp(\lambda(1 + x)) + \frac{\exp(2 \lambda)}{x^r}\right) C(y, z) F(x, y| z) g^n(s, y) g^n(s, z) dy dz ds dx \notag\\
& - \int_{\tfrac{1}{n}}^n \int_0^t   \left( \exp(\lambda(1 + x)) + \frac{\exp(2 \lambda)}{x^r}\right) g^n(s, x) \int_{\tfrac{1}{n}}^n y C(x, y) g^n(s, y) dy ds dx \notag\\
& \le \|g_0\|_{\lambda, r}^n + \int_{\tfrac{1}{n}}^n \int_0^t \int_{\tfrac{1}{n}}^n \int_x^n x  \left( \exp(\lambda(1 + x)) + \frac{\exp(2 \lambda)}{x^r}\right)\notag\\
&~~~ C(y, z) F(x, y| z) g^n(s, y) g^n(s, z) dy dz ds dx.
\end{align}
In the next, we focus on integration term of(\ref{to_show_f_n_in_Omega})
\begin{align}
& \int_{z = \tfrac{1}{n}}^n \int_{x = \tfrac{1}{n}}^n \int_{y = x}^n  \left( \exp(\lambda(1 + x)) + \frac{\exp(2 \lambda)}{x^r}\right) C(y, z) F(x, y| z) g^n(s, y) g^n(s, z) dy dz ds dx \notag\\
& = \int_{z = \tfrac{1}{n}}^n \int_{y = \tfrac{1}{n}}^n \int_{x = \tfrac{1}{n}}^y \left( \exp(\lambda(1 + x)) + \frac{\exp(2 \lambda)}{x^r}\right) C(y, z) F(x, y| z) g^n(s, y) g^n(s, z) dy dx \notag\\
& \le \left( \exp(\lambda(1 + n)) \lambda^{-1} + \frac{\exp(2 \lambda)n^{1 - r}}{1 - r}\right) \int_{z = \tfrac{1}{n}}^n \int_{y = \tfrac{1}{n}}^n   C(y, z) F(x, y| z) g^n(s, y) g^n(s, z) dy dx \notag\\
& \le \left( \exp(\lambda(1 + n)) \lambda^{-1} + \frac{\exp(2 \lambda)n^{1 - r}}{1 - r}\right) M^2 \int_{z = \tfrac{1}{n}}^n \int_{y = \tfrac{1}{n}}^n  y z g^n(s, y) g^n(s, z) dy dx \notag\\
& \le \left( \exp(\lambda(1 + n)) \lambda^{-1} + \frac{\exp(2 \lambda)n^{1 - r}}{1 - r}\right) M^2 (\|g_0\|_1)^2.
\end{align}
Therefore, from (\ref{to_show_f_n_in_Omega}), we obtain
\begin{align}
& \int_{\tfrac{1}{n}}^n \left( \exp(\lambda(1 + x)) + \frac{\exp(2 \lambda)}{x^r}\right) g^n(t, x) dx \notag\\
& \le \|g_0\|_{\lambda, r}^n + \left( \exp(\lambda(1 + n)) \lambda^{-1} + \frac{\exp(2 \lambda)n^{1 - r}}{1 - r}\right) M^2 (\|g_0\|_1)^2 T \notag\\
& = L.
\end{align}
As we have taken arbitrary $n$, hence the proof of Theorem \ref{thm_truncated} is accomplished.
\end{proof}
By the Theorem \ref{thm_truncated}, we draw the unique non-negative solution to (\ref{truncated_main_eq})-(\ref{truncated_ini_cond}) and denote it by $g^n, n = 1, 2, 3, ...$. For the complete domain, we consider the  ``zero" extension of each $g^n$, that is,

%
\begin{align}
 \hat{g}^n(t, x) =
  \begin{cases}
   g^n(t, x)   & \text{if } x\in[\tfrac{1}{n}, n] \\
   0 & \text{if } x \in(0, \infty) \setminus [\tfrac{1}{n}, n].
     \end{cases}
\end{align}
The $p^{\text{th}}$ order \textcolor{black}{truncated} moment of $\hat{g_n}(t, x)$ is denoted by
%
\begin{align}
N_{n, p}(t) = \int_0^\infty x^p \hat{g}^n(t, x) dx, p \in (- 1, \infty), n \in \mathbb{N}.
\end{align}
\begin{lemma}\label{boundedness_central_moment}
The moments are bounded \textcolor{black}{uniformly}, i.e.,
\begin{align}
N_{n, p}(t) = \int_0^\infty x^p \hat{g}^n(t, x) dx \le \bar{N}_p = \text{ constant if } 0 \le t \le T, ~ -1 < p < \infty \text{ and } n \geq 1.
\end{align}
%
\end{lemma}
\begin{proof}
For the first order moment, we integrate (\ref{reform_truncated_main_eq}) by multiplying $x$ and using the mass conservation law (\ref{truncated_mass_consv}), we get
\[\frac{d}{dt} N_{n, 1}(t) = 0.\]
So, we obtain
\begin{align}
N_{n, 1}(t) & = \int_0^\infty x \hat{g_0}(t, x) dx \le \int_0^\infty x g_0(x) dx = \textcolor{black}{\bar{N}_1, \text{a constant,}}~~ 0 \leq t \leq T, n \geq 1.\notag\\
\end{align}
All the integrations are exist as $C_n$ has compact support. By multiplying $x^2$ with (\ref{reform_truncated_main_eq}) and integrating, we obtain
\begin{align}
\frac{d N_{2, n}(t)}{dt} & = \int_0^\infty  \int_0^\infty \int_x^\infty x^2 C_n(y, z) F(x, y| z) g^n(t, y) g^n(t, z) dy dz dx \notag\\
                         & ~~~~~ ~~~~~ - \int_0^\infty \int_0^\infty x^2 g^n(t, x)  C(x, y) g^n(t, y) dy dx \notag\\
                         & \le \int_0^\infty  \int_0^\infty \left[\int_0^y x F(x, y| z) dx \right] y C_n(y, z) g^n(t, y) g^n(t, z) dy dz \notag\\
                         & ~~~~~ ~~~~~ - \int_0^\infty \int_0^\infty x^2 g^n(t, x)  C(x, y) g^n(t, y) dy dx \notag\\
                         & = \int_0^\infty  \int_0^\infty y^2 C_n(y, z) g_n(t, y) g^n(t, z) dy dz \notag\\
                         & ~~~~~ ~~~~~ - \int_0^\infty \int_0^\infty x^2 g^n(t, x)  C(x, y) g^n(t, y) dy dx \notag\\
                         & = 0.
\end{align}
Therefore, we get
\begin{align}
N_{2, n}(t) = \int_0^\infty x^2 \hat{g_0}(t, x) dx \le \int_0^\infty x^2 g_0(x) dx = \bar{N}_2 = \text{ constant, } 0 \leq t \leq T, n \geq 1.\notag\\
\end{align}
In the next, for $N_{n, -r}(t), ~ 0 < \sigma \le r < 1$, we note that
\begin{align}
\frac{d N_{n, -r}(t)}{dt} & = \int_0^n \int_{\tfrac{1}{n}}^n \int_x^n x^{-r} C_n(y, z) F(x, y| z) g^n(t, y) g^n(t, z) dy dz dx \notag\\                               & ~~ ~~~~~ - \int_0^n  \int_{\tfrac{1}{n}}^n x^{- r} g^n(t, x)  C(x, y) g^n(t, y) dy dx \notag\\
  & \le \int_0^n \int_{\tfrac{1}{n}}^n \int_x^n x^{- r} C_n(y, z) F(x, y| z) g^n(t, y) g^n(t, z) dy dz dx \notag\\
    & = \int_{z = \tfrac{1}{n}}^n \int_{y = 0}^n  \int_{x = 0}^y x^{- r} C_n(y, z) \frac{k_2}{y^\beta} g^n(t, y) g^n(t, z) dx dy dz \notag\\
     & = \int_{z = \tfrac{1}{n}}^n \int_{y = 0}^n \frac{y^{1 - r}}{1 - r} \frac{k_1 (1 + y)^\nu (1 + z)^\nu}{(y z)^\sigma}  \frac{k_2}{y^\beta} g^n(t, y) g^n(t, z) dx dy dz \notag\\
    & = \frac{k_1 k_2}{1 - r}  \int_{y = 0}^n y^{1 - r - \sigma - \beta} (1 + y)^\nu g^n(t, y) dy  \int_{z = \tfrac{1}{n}}^n (1 + z)^\nu z^{- \sigma}  g^n(t, z) dz \notag\\
   & \le \frac{k_1 k_2}{1 - r} \left( \int_{y = 0}^1 y^{1 - r - \sigma - \beta} (1 + y)^\nu g^n(t, y) dy + \int_{y = 1}^n y^{1 - r - \sigma - \beta} (1 + y)^\nu g^n(t, y) dy \right) \notag\\
    & \times \left(\int_{z = 0}^1 (1 + z)^\nu z^{- \sigma}  g^n(t, z) dz + \int_{z = 1}^n (1 + z)^\nu z^{- \sigma}  g^n(t, z) dz \right)\notag\\
    %
    & \le \frac{k_1 k_2}{1 - r} \left( \int_{y = 0}^1 (1 + y) y^{-r} g^n(t, y) dy + \int_{y = 1}^n y (1 + y)^\nu g^n(t, y) dy \right)  \notag\\
    & \times \left(\int_{z = 0}^1 (1 + z) z^{- r}  g^n(t, z) dz + \int_{z = 1}^n z (1 + z)  g^n(t, z) dz \right), \text{ as } 1-\sigma-\beta \ge 0, \notag\\
    %
    & \le \frac{k_1 k_2}{1 - r} ( 2 N_{-r, n}(t) + \bar{N}_1 + \bar{N}_2) (2 N_{-r, n}(t) + \bar{N}_1 + \bar{N}_2) \notag\\
    & =  \frac{k_1 k_2}{1 - r} (2 N_{-r, n}(t) + \bar{N}_1 + \bar{N}_2)^2.
\end{align}
Therefore,
\begin{align}
N_{n, -r}(t) = \int_0^\infty x^{-r} \hat{g}^n(t, x) dx \le \bar{N}_{-r} = \text{ constant if } 0 \le t \le T, ~ ~ 0 < \sigma \le r < 1 \text{ and } n \geq 1.
\end{align}
In order to achieve the boundedness of $N_{n, 0}(t)$, we observe that
\begin{align}\label{bound_N_0}
\frac{d N_{n, 0}(t)}{dt} & = \int_0^n \int_{\tfrac{1}{n}}^n \int_x^n C_n(y, z) F(x, y| z) g^n(t, y) g^n(t, z) dy dz dx \notag\\
                               & ~~~~~ ~~~~~ - \int_0^n  \int_{\tfrac{1}{n}}^n g^n(t, x)  C(x, y) g^n(t, y) dy dx \notag\\
                               & = \int_{z = \tfrac{1}{n}}^n \int_{y = 0}^n  \int_{x = 0}^y F(x, y| z) dx  C_n(y, z) g^n(t, y) g^n(t, z) dy dz \notag\\
                               & ~~~~~ ~~~~~ - \int_0^n  \int_{\tfrac{1}{n}}^n g^n(t, x)  C(x, y) g^n(t, y) dy dx \notag\\
                               & \le (N - 1) \int_0^n  \int_{\tfrac{1}{n}}^n g^n(t, x)  C(x, y) g^n(t, y) dy dx \notag\\
                               & \le (N - 1) \int_0^n  \int_{\tfrac{1}{n}}^n g^n(t, x) \frac{k_1 (1 + y)^\nu (1 + x)^\nu}{(y x)^\sigma}  g^n(t, y) dy dx \notag\\
                               & \le k_1  (N - 1) \int_{\tfrac{1}{n}}^n \frac{ (1 + x)^\nu}{ x^\sigma} g^n(t, x) dx  \int_0^n \frac{(1 + y)^\nu}{y^\sigma} g^n(t, y) dy.
\end{align}
For the last integration of (\ref{bound_N_0}), we observe that
\begin{align}
& \int_0^n \frac{(1 + y)^\nu}{y^\sigma} g^n(t, y) dy \notag\\
& =  \int_0^1 \frac{(1 + y)^\nu}{y^\sigma} g^n(t, y) dy +  \int_1^n \frac{(1 + y)^\nu}{y^\sigma} g^n(t, y) dy \notag\\
& \le 2^\nu \int_0^1 y^{- \sigma} g^n(t, y) dy +  \int_1^n (1 + y) g^n(t, y) dy \notag\\
&  \le 2^\nu \bar{N}_{- \sigma} + N_0(t) + \bar{N}_1.
\end{align}
Therefore, from (\ref{bound_N_0}), we obtain
\begin{align}
\frac{d N_{n, 0}(t)}{dt} \le k_1 (N - 1)(2^\nu \bar{N}_{- \sigma} + N_0(t) + \bar{N}_1)^2.
\end{align}
By integrating in the range $0$ to $t$, we obtain
\[N_{n, 0}(t) \le \bar{N}_0 = \text{ constant if } t \in [0, T] \text{ and } n \geq 1.\]
Now, for any number $q \in (-1, \infty)$, we note that
\begin{align}
N_{n, q}(t) & = \int_0^n x^q g^n(t, x) dx \notag\\
            & =   \int_0^1 x^q g^n(t, x) dx +  \int_1^n x^q g^n(t, x) dx \notag\\
            & \le N_{n, |q - 1|}(t) + N_{n, |q|}(t) \notag\\
            & \le \bar{N}_{n, |q - 1|} + \bar{N}_{n, |q|} = \bar{N}_q.
\end{align}
So, $N_{n, q}(t) \le \bar{N}_q$ (a constant independent of n), $t \in [0, T], n \ge 1$.\\
\par
\noindent Incorporating all above assessments, we get
\begin{align}\label{moments}
&~  N_{n, p}(t) \le \bar{N}_p = \text{a constant (independent of } t, n);~\text{for all}~  t \in [0, T],~n \in \mathbb{N} ~\text{ and }~ -1 < p < \infty.
\end{align}
\end{proof}
\begin{lemma}\label{relative_compact}
\textcolor{black}{In the topology functions (uniform convergence) constructed on each individual strip $\mathcal{P}(T, X_1, X_2)$, the sequence $\{\hat{g}^n\}_{n = 1}^\infty$ is relatively compact.}
\end{lemma}
\begin{proof}
Three steps are used to prove the lemma.\\
%
\textbf{Step 1.} At first, we prove that $\{\hat{g}^n\}_{n = 1}^\infty$ is \textcolor{black}{bounded uniformly} on $\mathcal{P}(T, X_1, X_2)$. By using (\ref{truncated_main_eq}), we have
\begin{align}\label{step_1_uniform_bounded}
\frac{\partial g^n(s, x)}{\partial s} & = \int_{\tfrac{1}{n}}^n \int_x^n C_n(y, z) F(x, y| z) g^n(s, y) g^n(s, z) dy dz \notag\\
                                      & ~~~~~ -  \int_{\tfrac{1}{n}}^n g^n(s, x)  C(x, y) g^n(s, y) dy \notag\\
                                      & \le \int_{\tfrac{1}{n}}^n \int_x^n C_n(y, z) F(x, y| z) g^n(s, y) g^n(s, z) dy dz \notag\\
                                      & \le k_1 k_2 \int_0^n \int_x^n \frac{(1 + y)^\nu (1 + z)^\nu}{(yz)^\sigma} \frac{1}{y^\beta} g^n(s, y) g^n(s, z) dy dz \notag\\
                                      &\le k_1 k_2 \frac{(1 + X_2)^{2 \nu}}{X_1^{\beta + 2 \sigma}} \int_0^n \int_x^n g^n(s, y) g^n(s, z) dy dz \notag\\
                                      & = k_1 k_2 \frac{(1 + X_2)^{2 \nu}}{X_1^{\beta + 2 \sigma}} N_{0, n}^2(t) \notag\\
\Rightarrow \frac{\partial g^n(s, x)}{\partial s} & \le k_1 k_2 \frac{(1 + X_2)^{2 \nu}}{X_1^{\beta + 2 \sigma}} \bar{N}_0^2. \notag\\
\text{Integrating the last} & \text{ inequality in the range } [0, t], \text{ we get}  g^n(t, x) \le k_1 k_2 \frac{(1 + X_2)^{2 \nu}}{X_1^{\beta + 2 \sigma}} \bar{N}_0^2 T = M_1 \text{ (say). }
\end{align}
\textbf{Step 2.}
In this part we show the equicontinuty of $\{\hat{g}^n\}_{n = 1}^\infty$  w.r.t. $t$ in every rectangular strip $\mathcal{P}(T, X_1, X_2)$.
\textcolor{black}{Fixing, $X_1 \le x \le X_2, 0 \le t \le t' \le T$} and recalling equation (\ref{truncated_main_eq}), we obtain
\begin{align}\label{step_2_equicontinuity_t}
|g_n(t', x) - g_n(t, x)| & \le \int_t^{t'} \left[\int_{\tfrac{1}{n}}^n \int_x^n C_n(y, z) F(x, y| z) g^n(s, y) g^n(s, z) dy dz \right.\notag\\
                                      & ~~~~~ \left. +  \int_{\tfrac{1}{n}}^n g^n(s, x)  C(x, y) g^n(s, y) dy \right] ds.
\end{align}
In the rectangle $\mathcal{P}(T, X_1, X_2)$, using boundedness of the moments (\ref{moments}), we get
\begin{align}
& \int_{\tfrac{1}{n}}^n \int_x^n C_n(y, z) F(x, y| z) g^n(s, y) g^n(s, z) dy dz \notag\\
& \le k_1 k_2 \int_{\tfrac{1}{n}}^n \int_x^n \frac{(1 + y)^\nu (1 + z)^\nu}{(yz)^\sigma} \frac{1}{y^\beta}  g^n(s, y) g^n(s, z) dy dz \notag\\
& \le k_1 k_2 \frac{(1 + X_2)^{2 \nu}}{X_1^{\beta + 2 \sigma}} \int_{\tfrac{1}{n}}^n \int_x^n g^n(s, y) g^n(s, z) dy dz \notag\\
& \le k_1 k_2 \frac{(1 + X_2)^{2 \nu}}{X_1^{\beta + 2 \sigma}} \bar{N}_0^2.
\end{align}
Therefore, the first integral becomes
\begin{align}
\int_t^{t'} \int_{\tfrac{1}{n}}^n \int_x^n C_n(y, z) F(x, y| z) g^n(s, y) g^n(s, z) dy dz ds
\le k_1 k_2 \frac{(1 + X_2)^{2 \nu}}{X_1^{\beta + 2 \sigma}} \bar{N}_0^2 (t' - t).
\end{align}
In the last integral of (\ref{step_2_equicontinuity_t}), applying uniform boundedness results (\ref{moments}) and (\ref{step_1_uniform_bounded}) we attain,
\begin{align}
& \int_t^{t'} \int_{\tfrac{1}{n}}^n g^n(s, x)  C(x, y) g^n(s, y) dy ds \notag\\
& \le k_1 M_1 \int_t^{t'} \int_{\tfrac{1}{n}}^n \frac{k_1 (1 + y)^\nu (1 + x)^\nu}{(y x)^\sigma} g^n(s, y) dy ds \notag\\
& \le k_1 M_1 \frac{(1 + X_2)^{2 \nu}}{X_1^{2 \sigma}} \bar{N}_0 (t' - t).
\end{align}
Therefore,
\begin{align}\label{eq_t_M2_(43)}
|g_n(t', x) - g_n(t, x)| \le M_2 (t' - t), \text{ where } M_2 = \left[ k_1 k_2 \frac{(1 + X_2)^{2 \nu}}{X_1^{\beta + 2 \sigma}} \bar{N}_0^2 + k_1 M_1 \frac{(1 + X_2)^{2 \nu}}{X_1^{2 \sigma}} \bar{N}_0\right].
\end{align}
We observe that $M_2$ is not dependent on \textcolor{black}{$t$ and $n$}. Hence, for an arbitrary $\epsilon > 0$, \textcolor{black}{we can choose $\delta (\le \frac{\epsilon}{M_2}) > 0$, so that}
\[ |g^n(t, x) - g^n(t', x)| < \epsilon \text{ whenever } |t - t'| < \delta, n \ge 1.\]
\end{proof}
\textbf{Step 3.}
At last step, we demonstrate the equicontinuity of $\{\hat{g}^n\}_{n = 1}^\infty$  with respect to $x$ in each rectangle $\mathcal{P}(T, X_1, X_2)$. Let us suppose, $X_2 \ge x' \ge x \ge X_1$, then from (\ref{truncated_main_eq}) for every $n \geq 1$, we have
\begin{align}\label{equi_cont_x_main_eq}
|g^n(t, x') - g^n(t, x)| & \le |g_0(x') - g_0(x)| \notag\\
                         & + \int_0^t \int_{\tfrac{1}{n}}^n |C_n(x', y) - C_n(x, y)| g^n(s, x') g^n(s, y) dy ds\notag\\
                         & + \int_0^t \int_{\tfrac{1}{n}}^n  C_n(x, y) |g^n(s, x') -g^n(s, x)| g^n(s, y) dy ds\notag\\
                         & + \int_0^t \int_{\tfrac{1}{n}}^n\int_{x'}^n C_n(y, z) |b(x', y; z) - F(x, y| z)| g^n(s, y) g^n(s, z) dy dz ds\notag\\
                         & + \int_0^t \int_{\tfrac{1}{n}}^n \int_{x'}^x C_n(y, z) F(x, y| z) g^n(s, y) g^n(s, z) dy dz ds.
\end{align}
From the first integral of (\ref{equi_cont_x_main_eq}), we observe that
\begin{align}\label{step_3_(51)}
& \int_0^t \int_{\tfrac{1}{n}}^n |C_n(x', y) - C_n(x, y)| g^n(s, x') g^n(s, y) dy ds\notag\\
& \le \underbrace{ \int_0^t \int_0^{z_1} |C_n(x', y) - C_n(x, y)| g^n(s, x') g^n(s, y) dy ds}_{= I_1}\notag\\
& + \underbrace{ \int_0^t \int_{z_1}^{z_2} |C_n(x', y) - C_n(x, y)| g^n(s, x') g^n(s, y) dy ds}_{= I_2}\notag\\
& +  \underbrace{ \int_0^t \int_{z_2}^\infty |C_n(x', y) - C_n(x, y)| g^n(s, x') g^n(s, y) dy ds}_{= I_3}.
\end{align}
Therefore, from (\ref{step_3_(51)}), we obtain
\begin{align}
I_1 & = \int_0^t \int_0^{z_1} |C_n(x', y) - C_n(x, y)| g^n(s, x') g^n(s, y) dy ds\notag\\
    & = k_1 \int_0^t \int_0^{z_1} \left|\frac{ (1 + x')^\nu (1 + y)^\nu }{(x' y )^\sigma} - \frac{  (1 + y)^\nu (1 + x)^\nu}{(y x)^\sigma}  \right|g^n(s, x') g^n(s, y) dy ds \notag\\
    & = k_1 \int_0^t \int_0^{z_1} \frac{(1 + y)^\nu}{y^\sigma} \left|\frac{ (1 + x')^\nu}{(x')^\sigma} - \frac{(1 + x)^\nu}{x^\sigma}  \right|g^n(s, x') g^n(s, y) dy ds \notag\\
    & \le k_1 (1 + z_1)^\nu \frac{(1 + X_2)^\nu}{X_1^\sigma} \int_0^t \int_0^{z_1} y^{ -\sigma} g^n(s, x') g^n(s, y) dy ds \notag\\
    & \le k_1 (1 + z_1)^\nu \frac{(1 + X_2)^\nu}{X_1^\sigma} M_1 T \int_0^{z_1} y^{ -\sigma} g^n(s, y) dy.
\end{align}
We choose $z_1$ such that $\int_0^{z_1} y^{ -\sigma} g^n(s, y) dy < \epsilon $, it is possible because $\int_0^\infty y^{ -\sigma} g^n(s, y) dy < \bar{N}_{-\sigma} < \infty$ and therefore, $\int_0^{z_1} y^{ -\sigma} g^n(s, y) dy \rightarrow 0 $ as $z_1 \rightarrow 0.$
In the next, from (\ref{step_3_(51)}), we obtain
\begin{align}
I_2 & = \int_0^t \int_{z_1}^{z_2} |C_n(x', y) - C_n(x, y)| g^n(s, x') g^n(s, y) dy ds\notag\\
    & = \epsilon \int_0^t \int_{z_1}^{z_2} g^n(s, x') g^n(s, y) dy ds \notag\\
    & \le \epsilon M_1 T \bar{N}_0.
\end{align}
Finally, from (\ref{step_3_(51)}), we get
\begin{align}
I_3 & = \int_0^t \int_{z_2}^\infty |C_n(x', y) - C_n(x, y)| g^n(s, x') g^n(s, y) dy ds\notag\\
    & = k_1 \int_0^t \int_{z_2}^\infty \left|\frac{ (1 + x')^\nu (1 + y)^\nu }{(x' y )^\sigma} - \frac{  (1 + y)^\nu (1 + x)^\nu}{(y x)^\sigma}  \right|g^n(s, x') g^n(s, y) dy ds \notag\\
    & = k_1 \int_0^t \int_{z_2}^\infty \frac{(1 + y)^\nu}{y^\sigma} \left|\frac{ (1 + x')^\nu}{(x')^\sigma} - \frac{(1 + x)^\nu}{x^\sigma}  \right|g^n(s, x') g^n(s, y) dy ds \notag\\
    & \le k_1 \frac{(1 + X_2)^\nu}{X_1^\sigma} \int_0^t \int_{z_2}^\infty \frac{y^{1 + \sigma}}{z_2^{1 + \sigma}} \frac{(1 + y)}{y^\sigma} g^n(s, x') g^n(s, y) dy ds \notag\\
    & \le k_1 \frac{(1 + X_2)^\nu}{X_1^\sigma} \frac{1}{z_2^{1 + \sigma}} M_1 T (\bar{N}_1 + \bar{N}_2).
\end{align}
Choose $z_2$ such that $\frac{(1 + X_2)^\nu}{X_1^\sigma} (\bar{N}_1 + \bar{N}_2) \frac{1}{z_2^{1 + \sigma}} < \epsilon.$
For the third integral of (\ref{equi_cont_x_main_eq}), we get
\begin{align}\label{J_1_J_2}
& \int_0^t \int_{\tfrac{1}{n}}^n \int_{x'}^n C_n(y, z) |b(x', y; z) - F(x, y| z)| g^n(s, y) g^n(s, z) dy dz ds \notag\\
& = \underbrace{ \int_0^t \int_{\tfrac{1}{n}}^n \int_{x'}^{z_2} C_n(y, z) |b(x', y; z) - F(x, y| z)| g^n(s, y) g^n(s, z) dy dz ds}_{J_1} \notag\\
& ~~~~ + \underbrace{ \int_0^t \int_{\tfrac{1}{n}}^n \int_{z_2}^n C_n(y, z) |b(x', y; z) - F(x, y| z)| g^n(s, y) g^n(s, z) dy dz ds }_{J_2}.
\end{align}
From (\ref{J_1_J_2}), we obtain following inequality
\begin{align}
J_1 & \le \epsilon \int_0^t \int_{\tfrac{1}{n}}^n \int_{x'}^{z_2} C_n(y, z) g^n(s, y) g^n(s, z) dy dz ds \notag\\
    & \le \epsilon \int_0^t \int_{\tfrac{1}{n}}^n \int_{x'}^{z_2} k_1 \frac{(1 + y)^\nu (1 + z)^\nu}{(y z)^\sigma} g^n(s, y) g^n(s, z) dy dz ds \notag\\
    & \le k_1 \epsilon \int_0^t \int_{\tfrac{1}{n}}^n \int_{x'}^{z_2} (y^{- \sigma} + y^{\nu - \sigma}) g^n(s, y) (z^{- \sigma} + z^{\nu - \sigma}) g^n(s, z) dy dz ds \notag\\
    & \le k_1 \epsilon (\bar{N}_{- \sigma} + \bar{N}_{\nu - \sigma})^2 T.
\end{align}
From (\ref{J_1_J_2}), we obtain following inequality
\begin{align}
J_2 & = \int_0^t \int_{\tfrac{1}{n}}^n \int_{z_2}^n C_n(y, z) |b(x', y; z) - F(x, y| z)| g^n(s, y) g^n(s, z) dy dz ds \notag\\
    & = k_1 k_2 \int_0^t \int_{\tfrac{1}{n}}^n \int_{z_2}^n \frac{(1 + y)^\nu (1 + z)^\nu}{(y z)^\sigma} \frac{1}{y^\beta} g^n(s, y) g^n(s, z) dy dz ds \notag\\
    & \le  k_1 k_2 z_2^{- \beta} \int_0^t \int_{\tfrac{1}{n}}^n \int_{z_2}^n (y^{- \sigma} + y^{\nu - \sigma}) g^n(s, y) (z^{- \sigma} + z^{\nu - \sigma}) g^n(s, z) dy dz ds \notag\\
    & =  k_1 k_2 z_2^{- \beta} (\bar{N}_{- \sigma} + \bar{N}_{\nu - \sigma})^2 T.
\end{align}
Choose $z_2$ such that $z_2^{- \beta} (\bar{N}_{- \sigma} + \bar{N}_{\nu - \sigma})^2 < \epsilon.$\\
The second integral of (\ref{equi_cont_x_main_eq}) reduces to
\begin{align}
& \int_0^t \int_0^n C_n(x, y) |g^n(s, x') - g^n(s, x)| g^n(s, y) dy ds \notag\\
& \le k_1 \int_0^t \int_0^n \frac{(1 + x)^\nu (1 + y)^\nu}{(x y)^\sigma} w_n(s) g^n(s, y) dy ds \notag\\
& \le k_1 \frac{(1 + X_2)^\nu}{X_1^\sigma} \int_0^t \int_0^n 2^\nu (y^{- \sigma} + y^{\nu - \sigma}) w_n(s) g^n(s, y) dy ds \notag\\
& \le 2^\nu k_1 \frac{(1 + X_2)^\nu}{X_1^\sigma} (\bar{N}_{- \sigma} + \bar{N}_{\nu - \sigma})\int_0^t w_n(s) ds. \notag\\
\end{align}
Now, we estimate last integral of (\ref{equi_cont_x_main_eq}) as
\begin{align}
& \int_0^t \int_0^n \int_x^{x'} C_n(y, z) F(x, y| z) g^n(s, y) g^n(s, z) dy dz ds \notag\\
& = \int_0^t \int_0^n \int_0^{x'} C_n(y, z) F(x, y| z) g^n(s, y) g^n(s, z) dy dz ds \notag\\
& ~~~~~ - \int_0^t \int_0^n \int_0^x C_n(y, z) F(x, y| z) g^n(s, y) g^n(s, z) dy dz ds.
\end{align}
The last two integrals in $\mathcal{P}(T, X_1, X_2)$ are continuous due to compact support.
Therefore, $\exists$ some $M_6 > 0$, so that $|x - x'| < \delta$ indicates
\begin{align}
\int_0^t \int_0^n C_n(x, y) |g^n(s, x') - g^n(s, x)| g^n(s, y) dy ds < M_6 \epsilon.
\end{align}
Taking all of the result into account, we arrive at
%
%
\[ w_n (t) \le M_3 \epsilon + M_4 \epsilon + M_5 \int_0^t w_n(s) ds + M_6 \epsilon. \]
All the above $M_i$'s $(i=3,4,5,6)$ are independent of $t$ and $n$, and using Gronwall's inequality, we get
\begin{align}\label{w_n_(59)}
w_n(t) \le (M_3 + M_4 + M_6) \epsilon \exp(M_5 T) = M_7 \epsilon,
\end{align}
where $M_7$ is a positive constant. Hence, the claim of Step 3 is proved. We conclude from (\ref{eq_t_M2_(43)}) and (\ref{w_n_(59)}) that
\begin{align}
\sup_{|x' - x| < \delta, |t' - t|< \delta} |g^n(t', x') - g^n(t, x)| \le (M_2 + M_7) \epsilon, X_1 \le x, x'\le X_2, 0 \le t, t' < T.
\end{align}
\textcolor{black}{Hence by Arzel\`{a}-Ascoli theorem, Lemma \ref{relative_compact} follows.}

\begin{theorem}\textbf{(Global in-time existence of a solution)}\label{thm 1}\\
Assume the symmetric collision kernel $C(x, y)$ and the fragmentation rate $F(x, y| z)$ be \textcolor{black}{non-negative and continuous} in $(0, \infty) \times (0, \infty)$  and $(0, \infty) \times (0, \infty) \times (0, \infty)$ respectively. Moreover,
we assume that
\begin{enumerate}[$($i$)$]
\item $C(x,y) \le k_1~ \frac{(1 + x)^\nu (1 + y)^\nu}{\left(xy\right)^\sigma}$, for all ~$x, y \in (0,\infty),$ where $k_1 > 0$ is a constant, $\sigma \in \left[0,\tfrac{1}{2}\right]$ and $\mu \in [0, 1]$,
\item for all $0 < x < y,$ there exist real number $0 < \beta \le \sigma$, so that $F(x, y| z) \le \frac{k_2}{y^\beta}$, where $k_2$ is a positive constant.
\end{enumerate}
Let the initial data function obey $g_0 \in \Omega_{.,r_2}^+ (0)$. Then there exists a solution in $\Omega_{.,r_2}^+ (T)$  of the problem (\ref{main_eq})-(\ref{ini_cond}). \\
\end{theorem}
\begin{proof}
Using diagonal method \cite{paul2018existence}, there exist a subsequence $\{\hat{g}^i\}_{i = 1}^\infty$ from $\{\hat{g}^n\}_{n = 1}^\infty$ converging to a \textcolor{black}{non-negative and continuous} function $g$ that obeys (\ref{moments}), on each compactly supported region $\mathcal{P}$.\\
Now, we focus on the integral $\int_{x_1}^{x_2} x^p g(t,x) dx$, for $-1 < p < \infty.$
As, $\exists,~ i \geq 1$, so that for $\epsilon > 0, 0 < x_1 < x_2 < \infty, -1 < p < \infty,$
\begin{align}
\int_{x_1}^{x_2} x^p g(t, x) dx \le \int_{x_1}^{x_2} x^p \hat{g}^i(t, x) dx \epsilon \le \bar{N}_p + \epsilon;
\end{align}
therefore,
\begin{align}\label{moments_f}
\int_0^\infty x^p g(t, x) dx \le \bar{N}_p, p \in (-1, \infty),
\end{align}
as $\epsilon, x_1, x_2$ are arbitrary. In the similar way,
\[ \int_0^\infty x g(t, x) dx \le \bar{N}_1 = \int_0^\infty x g_0(x) dx. \]
Finally, we aim the limiting function $g(t, x)$ is a solution of (\ref{main_eq}) and (\ref{ini_cond}). By replacing $C_i$ and $\hat{g}^i$ in (\ref{truncated_main_eq}) with $C_i - C + C$ and $\hat{g}^i - g + g$, respectively
\begin{align}\label{conv_thm_2(64)}
&(\hat{g}^i - g)(t, x) + g(t, x)  = g_0(x) + \int_0^t \bigg[ \int_0^\infty \int_x^\infty (C_i - C)(y, z) F(x, y| z) \hat{g}^i(s, y) \hat{g}^i(s, z) dy dz \notag\\
                                & + \int_0^\infty \int_x^\infty C(y, z) F(x, y| z) (\hat{g}^i - g)(s, y) \hat{g}^i(s, z) dy dz \notag\\
                                & + \int_0^\infty \int_x^\infty C(y, z) F(x, y| z) g(s, y) (\hat{g}^i - g)(s, z) dy dz - \hat{g}^i(s, x) \int_0^\infty (C_i - C)(x, y) \hat{g}^i(s, y) dy \notag\\
                                & - (\hat{g}^i - f)(s, x) \int_0^\infty C(x, y) \hat{g}^i(s, y) dy -g(s, x) \int_0^\infty C(x, y) (\hat{g}^i - g)(s, y) dy \notag\\
                                & + \int_0^\infty \int_x^\infty C(y, z) F(x, y| z) g(s, y) g(s, z) dy dz - g(s, x) \int_0^\infty C(x, y) g(s, y) dy
                                \bigg] ds.
\end{align}
In the next, passing the limit $i \rightarrow \infty$ in (\ref{conv_thm_2(64)}), We can see that, with the exception of the last two, all of the infinite integral terms tend to zero due to their estimates. To demonstrate this, we'll do the following steps:
\begin{align}\label{3.76}
& \left| \int_0^\infty \int_x^\infty (C_i - C)(y, z) F(x, y| z) \hat{g}^i(s, y) \hat{g}^i(s, z) dy dz \right| \notag\\
& \le \left| \int_0^\infty \int_x^{y_1} (C_i - C)(y, z) F(x, y| z) \hat{g}^i(s, y) \hat{g}^i(s, z) dy dz \right| \notag\\
& ~~~~ + \left| \int_0^\infty \int_{y_1}^{y_2} (C_i - C)(y, z) F(x, y| z) \hat{g}^i(s, y) \hat{g}^i(s, z) dy dz \right| \notag\\
& ~~~~ + \left| \int_0^\infty \int_{y_2}^\infty (C_i - C)(y, z) F(x, y| z) \hat{g}^i(s, y) \hat{g}^i(s, z) dy dz \right|.
\end{align}
The first estimate in the r.h.s of \eqref{3.76} is calculated as
\begin{align}
& \left| \int_0^\infty \int_x^{y_1} (C_i - C)(y, z) F(x, y| z) \hat{g}^i(s, y) \hat{g}^i(s, z) dy dz \right| \notag\\
& \le \left| \int_0^{z_1} \int_x^{y_1} (C_i - C)(y, z) F(x, y| z) \hat{g}^i(s, y) \hat{g}^i(s, z) dy dz \right| \notag\\
& + \left| \int_{z_1}^{z_2} \int_x^{y_1} (C_i - C)(y, z) F(x, y| z) \hat{g}^i(s, y) \hat{g}^i(s, z) dy dz \right| \notag\\
& + \left| \int_{z_2}^\infty \int_x^{y_1} (C_i - C)(y, z) F(x, y| z) \hat{g}^i(s, y) \hat{g}^i(s, z) dy dz \right|.
\end{align}
In the similar fashion as (\ref{step_3_(51)}), the integral shown above can be made as small as desirable. Similarly, the remaining two integrals of (\ref{conv_thm_2(64)}) are equally small. Hence, $\exists$ a $p_1>0$ so that for all $i \geq p_1$,
\begin{align}\label{conv_thm_2(67)}
& \left| \int_0^\infty \int_x^\infty (C_i - C)(y, z) F(x, y| z) \hat{g}^i(s, y) \hat{g}^i(s, z) dy dz \right| < M_8 \epsilon.
\end{align}
Applying the same logic, we can find two numbers $p_k(k=1,2)>0$, such that
\begin{align}\label{conv_thm_2(68-69)}
&  \int_0^\infty \int_x^\infty C(y, z) F(x, y| z) (\hat{g}^i - g)(s, y) \hat{g}^i(s, z) dy dz \le M_9 \epsilon,   \forall~ i \geq p_2, \notag\\
&  \int_0^\infty \int_x^\infty C(y, z) F(x, y| z) g(s, y) (\hat{g}^i - g)(s, z) dy dz \le M_{10}  \epsilon, \forall~ i \geq p_3.
\end{align}
Taking $p = \max\{p_1, p_2, p_3\}$, all the preceding integrals tends to zero as $\epsilon$ can be made sufficiently small. Continuing the above procedure, one can easily demonstrate that other expressions in (\ref{conv_thm_2(64)}) tends to zero. Hence, the limiting solution $g(t, x)$ satisfies (\ref{main_eq}) and (\ref{ini_cond}), written in integral form:
\begin{align}\label{conclusion_eq_thm_2_(70)}
g(t, x)  & =  g_0(x) + \int_0^t \bigg[\int_0^\infty \int_x^\infty C(y, z) F(x, y| z) g(s, y) g(s, z) dy dz \notag\\
         &                              - g(s, x) \int_0^\infty C(x, y) g(s, y) dy\bigg] ds.
\end{align}
Therefore, from (\ref{conv_thm_2(67)})-(\ref{conv_thm_2(68-69)}), and the continuity property of $g(t, x)$ indicates that the right side of (\ref{conclusion_eq_thm_2_(70)}) is continuous on $\mathcal{P} = \{(t, x): t \in [0, T], 0 < x < \infty\}$. Moreover, differentiating (\ref{conclusion_eq_thm_2_(70)}) respect to $t$ establishes $g(t, x)$ is a differentiable continuous solution of (\ref{main_eq} - \ref{ini_cond}) and the condition (\ref{moments_f}) concludes that $g(t, x)$ remains in $\Omega_{.,r}^+(T)$.
\end{proof}

\begin{theorem}\textbf{(Mass conservation of global solution)}\label{Mass_conv_thm}\\
The solution of (\ref{main_eq}) and (\ref{ini_cond}) satisfies the conservation of mass property in-line with the conditions of Theorem \ref{thm 1}.
\end{theorem}
\begin{proof}
We denote the mass $\mathfrak{M} = \int_0^\infty x g(t, x) dx.$
By multiplication (\ref{main_eq}) with an weight $x$ and then integrating it, we get
\begin{align}\label{eq_for_mass_conv}
\frac{d\mathfrak{M}}{dt} & =  \int_0^\infty x \int_0^\infty \int_x^\infty C(y, z) F(x, y| z) g(t, y) g(t, z) dy dz dx \notag\\
              & - \int_0^\infty x g(t, x) \int_0^\infty C(x, y) g(t, y) dy dx \notag\\
              & = \mathfrak{M}_1 - \mathfrak{M_2} \text{ (say).}
\end{align}
A change in the order of integration of $\mathfrak{M}_1$ yields
\begin{align}
\mathfrak{M}_1 & =  \int_0^\infty x \int_0^\infty \int_x^\infty C(y, z) F(x, y| z) g(t, y) g(t, z) dy dz dx \notag\\
               & = \int_0^\infty \int_0^\infty C(y, z) g(t, y) g(t, z) dy dz \int_0^y x F(x, y| z) dx \notag\\
               & = \int_0^\infty  \int_0^\infty y C(y, z) g(t, y) g(t, z) dy dz   \notag\\
               & = \int_0^\infty \int_0^\infty x C(x, y) g(t, x) g(t, y) dy dx  = \mathfrak{M_2}.
\end{align}
In the next, we show that the integral of $\mathfrak{M_2}$ are finite. Here
\begin{align}\label{bdd_mass_conv}
& \int_0^\infty \int_0^\infty x C(x, y) g(t, x) g(t, y) dy dx \notag\\
& \le \int_0^\infty \int_0^\infty x k_1 \frac{(1 + x)^\nu (1 + y)^\nu} {(xy)^\sigma} g(t, x) g(t, y) dy dx \notag\\
& \le k_1 \int_0^\infty  (x^{1 - \sigma} + x^{1 + \nu -\sigma}) g(t, x) dx \int_0^\infty(y^{- \sigma} + y^{\nu - \sigma}) g(t, y) dy \notag\\
& \le k_1 (\bar{N}_{1 - \sigma} + \bar{N}_{1 + \nu -\sigma}) (\bar{N}_{- \sigma} + \bar{N}_{\nu -\sigma}) < \infty.
\end{align}
\textcolor{black}{Therefore, using Fubini's integral theorem, the two integrals $\mathfrak{M_1}$ and $\mathfrak{M_2}$ of (\ref{eq_for_mass_conv}) are finite quantity in nature, and thereafter, $$\frac{d\mathfrak{M}}{dt} = 0.$$} Hence, the mass conservation law holds.

\end{proof}

\section{Uniqueness of global solution}
\begin{theorem}\label{thm_uniqueness}
Let the symmetric collision kernel $K(x, y)$ and the fragmentation rate $b(x, y; z)$ be non-negative and continuous in $(0, \infty) \times (0, \infty)$  and $(0, \infty) \times (0, \infty) \times (0, \infty)$ respectively. Moreover,
we assume that
\begin{enumerate}[$($i$)$]
\item for all ~$x, y \in (0,\infty), ~C(x,y) \le k_1~ \frac{1}{\left(xy\right)^\sigma}$, where $k_1 > 0$ is a constant, $\sigma \in \left[0,\tfrac{1}{2}\right]$,
\item for all $y > x > 0,$ there exist real number $0 < \beta \le \sigma$, so that $b(x, y; z) \le \frac{k_2}{y^\beta}$, where $k_2$ is a positive constant.
\end{enumerate}
In addition, let the initial data function satisfy $f_0 \in \Omega_{.,r}^+ (0)$. Then the solution to the problem (\ref{main_eq}) with (\ref{ini_cond}) is unique in $\Omega_{.,r} (T)$. \\
\end{theorem}
\begin{proof}
If possible, let there exist two mass-conserving solutions $g_1(t, x)$ and $g_2(t, x)$ in $\Omega_{.,r}(T)$ to the IVP (\ref{main_eq}) with (\ref{ini_cond}). We will show that $g_1 = g_2$. Let $\phi(t, x) = g_1(t, x) -g_2(t, x)$ and $\psi(t, x) = g_1(t, x) + g_2(t, x)$. Since $g_1(t, x), g_2(t, x) \in \Omega_{.,r}(T)$, there exists $\hat{\lambda} > 0$ such that
\begin{align}\label{boundedness_g_i}
\int_0^{\infty} \left( \exp(\hat{\lambda }( 1 + x) ) + \frac{\exp(2 \hat{\lambda })}{x^{r}}\right)|g_i(t, x)| dx  ~< + \infty.
\end{align}
uniformly w.r.t. $t \in [0,T]$. Let $0 \le \lambda < \hat{\lambda}$. Then from the definition of $\phi$, we have
\begin{align}\label{Uniq:eq_small_phi}
    \frac{\partial \phi(t, x)}{\partial t} = & \int_0^\infty \int_x^\infty C(y, z) F(x,y|z) [g_1(t, y) g_1(t, z) - g_2(t, y) g_2(t, z)] dy dz \notag\\
    & - \int_0^\infty C(x, y) [g_1(t, x)g_1(t, y) - g_2(t, x) g_2(t,y)] dy.
\end{align}
We define
\begin{align}\label{defn_phi_psi}
    \Phi(t, \lambda) & = \int_0^\infty \left[ \exp(\lambda x) + \frac{1}{x^\theta}\right]  \phi(t, x) dx\notag\\
    \text{ and } \Psi(t, \lambda) & = \int_0^\infty \left[ \exp(\lambda x) + \frac{1}{x^\theta}\right]  \psi(t, x) dx,
\end{align}
where $\theta + \sigma <1.$
From (\ref{Uniq:eq_small_phi}), $\Phi(t, \lambda)$ obey
\begin{align}\label{Uq_eq_Phi}
    \Phi(t, \lambda)  = & \int_0^t \int_{x = 0}^\infty \left[ \exp(\lambda x) + \frac{1}{x^\theta}\right] \sgn(\phi(s, x))\notag\\
    & \left[\int_{z = 0}^\infty \int_{y = x}^\infty C(y, z) F(x,y|z) [g_1(s, y) g_1(s, z) - g_2(s, y) g_2(s, z)] dy dz\right. \notag\\
    & \left.- \int_{y = 0}^\infty C(x, y) [g_1(s, x)g_1(s, y) - g_2(s, x) g_2(s,y)] dy\right] dx ds.
\end{align}
Let,
\begin{align}
    G(t, y, z) = & [g_1(t, y) g_1(t, z) - g_2(t, y) g_2(t, z)] \notag\\
    = & [g_1(t, y) - g_2(t, y)] g_1(t, z) + g_2(t, y) [g_1(t, z) - g_2(t, z)] \notag\\
    = & \phi(t,y) g_1(t, z) + g_2(t,y) \phi(t, z).
\end{align}
From (\ref{Uq_eq_Phi}), we obtain
\begin{align}\label{uniq_break}
    & \int_{x = 0}^\infty \left[ \exp(\lambda x) + \frac{1}{x^\theta}\right]\sgn(u(s, x))  \int_{z = 0}^\infty \int_{y = x}^\infty C(y, z) F(x,y|z) G(t, y, z) dy dz \notag\\
    & = \int_{z = 0}^\infty\int_{y = 0}^\infty \int_{x = 0}^y \left[ \exp(\lambda x) + \frac{1}{x^\theta}\right] k_1 (yz)^{- \sigma} \frac{k_2}{y^\beta} |G(t, y, z)| dx dy dz\notag\\
    & \le k_1 k_2  \int_{z = 0}^\infty\int_{y = 0}^\infty \left[ y \exp(\lambda y) + \frac{y^{1-\theta}}{1 - \theta}\right] y^{-\sigma - \beta} z^{-\sigma}|G(t, y, z)|  dy dz\notag\\
    & \le k_3 \int_{z = 0}^\infty\int_{y = 0}^\infty \left[  \exp(\lambda y) + y^{-\theta}\right] y^{1-\sigma - \beta} z^{-\sigma}|G(t, y, z)|  dy dz,
\end{align}
where $k_3 = k_1 k_2 \max\{1, (1-\theta)^{-1}\}$.
Now from (\ref{uniq_break}) we get
\begin{align}
    & k_3 \int_{z = 0}^\infty\int_{y = 0}^\infty \left[  \exp(\lambda y) + y^{-\theta}\right] y^{1-\sigma - \beta} z^{-\sigma}|G(t, y, z)|  dy dz \notag\\
    & = k_3 \int_{z = 0}^\infty\int_{y = 0}^1 \left[  \exp(\lambda y) + y^{-\theta}\right] y^{1-\sigma - \beta} z^{-\sigma}|G(t, y, z)|  dy dz \notag\\
    & + k_3 \int_{z = 0}^\infty\int_{y = 1}^\infty \left[  \exp(\lambda y) + y^{-\theta}\right] y^{1-\sigma - \beta} z^{-\sigma}|G(t, y, z)|  dy dz \notag\\
    & \le k_3 \int_{z = 0}^\infty\int_{y = 0}^1 \left[  \exp(\lambda y) + y^{-\theta}\right]  z^{-\sigma}|G(t, y, z)|  dy dz \notag\\
    & + k_3 \int_{z = 0}^\infty\int_{y = 1}^\infty  [\exp(\lambda y) +1] y z^{-\sigma}|G(t, y, z)|  dy dz \notag\\
    & =  k_3 \int_{z = 0}^\infty\int_{y = 0}^1 \left[  \exp(\lambda y) + y^{-\theta}\right]  z^{-\sigma}|\phi(t,y) g_1(t, z) + g_2(t,y) \phi(t, z)|  dy dz \notag\\
    &+ k_3 \int_{z = 0}^\infty\int_{y = 1}^\infty  2\exp(\lambda y)  y z^{-\sigma}|\phi(t,y) g_1(t, z) + g_2(t,y) \phi(t, z)|  dy dz \notag\\
    & \le k_3 [\Phi(t, \lambda,) \Psi(t, \lambda)+ \Psi(t, \lambda)\Phi(t, \lambda)] + 2k_3 [\Phi_\lambda(t, \lambda) \Psi(t, \lambda) + \Psi_\lambda(t, \lambda) \Phi(t, \lambda)].
\end{align}
From (\ref{boundedness_g_i}) we can conclude that
\begin{align}
    \Psi(t, \lambda)=  \int_0^\infty \left[ \exp(\lambda x) + \frac{1}{x^\theta}\right]  (g_1 (t, x) + g_2(t, x)) dx < \bar{B}_1, \text{ a constant.}
\end{align}
We choose small $\epsilon_1 >0$ such that $\lambda + \epsilon_1 < \hat{\lambda}$, therefore
\begin{align}
    \Psi_\lambda(t, \lambda) & = \int_{x = 0}^\infty x \exp(\lambda x) (g_1(t, x) + g_2(t, x)) dx \notag\\
    & \le \int_{x = 0}^\infty  \exp((\lambda + \epsilon_1)x) (g_1(t, x) + g_2(t, x)) dx \notag\\
    &\le \int_{x = 0}^\infty  \exp(\hat{\lambda}x) (g_1(t, x) + g_2(t, x)) dx \notag\\
    & = \bar{B}_2 ,\text{ a constant} < \infty.
\end{align}
For the last integral of (\ref{Uq_eq_Phi}), we have
\begin{align}
    &- \int_{x = 0}^\infty\int_{y = 0}^\infty\left[\exp(\lambda x) + \frac{1}{x^\theta}\right] \sgn(\phi(s, x)) C(x, y) [g_1(s, x)g_1(s, y) - g_2(s, x) g_2(s,y)] dy dx \notag\\
    =&  - \int_{x = 0}^\infty\int_{y = 0}^\infty\left[\exp(\lambda x) + \frac{1}{x^\theta}\right]\sgn(\phi(s, x)) C(x, y) [\phi(s, x)g_1(s, y) + g_2(s, x) \phi(s,y)] dy dx\notag\\
    = &  \underbrace{- \int_{x = 0}^\infty\int_{y = 0}^\infty\left[\exp(\lambda x) + \frac{1}{x^\theta}\right] C(x, y) |\phi(s, x)|g_1(s, y) dy dx }_{\le 0} \notag\\
    &- \int_{x = 0}^\infty\int_{y = 0}^\infty\left[\exp(\lambda x) + \frac{1}{x^\theta}\right]C(x, y)g_2(s, x) \sgn(\phi(s, x))\phi(s,y) dy dx\notag\\
    \le & k_1 \int_{x = 0}^\infty\int_{y = 0}^\infty\left[\exp(\lambda x) + \frac{1}{x^\theta}\right] x^{-\sigma} y^{- \sigma}g_2(s, x) |\phi(s,y)| dy dx\notag\\
    \le & k_1 \left[\int_{x = 0}^1\left[\exp(\lambda x) + \frac{1}{x^\theta}\right] x^{-\sigma} g_2(s, x)   dx+ \int_{x =1 }^\infty\left[\exp(\lambda x) + \frac{1}{x^\theta}\right] x^{-\sigma} g_2(s, x)  dx\right]\notag\\
     & \int_{y = 0}^\infty y^{- \sigma}|\phi(s,y)|dy\notag\\
    \le & k_1 \left[\int_{x = 0}^1\left[x^{-\sigma} + x^{-\theta -\sigma}\right]  g_2(s, x)   dx+ \int_{x =1 }^\infty\left[\exp(\lambda x) + \frac{1}{x^\theta}\right]  g_2(s, x)  dx\right]\Phi(\lambda, s )\notag\\
    \le & k_1 [\bar{N}_{-\sigma} + \bar{N}_{-\sigma - \theta} + \Psi(s, \lambda)]\Phi(\lambda, s ).
\end{align}
From (\ref{Uq_eq_Phi}), we get
\begin{align}\label{ineq_phi}
    \Phi(t, \lambda,) & \le \int_0^t 2 k_3 [\Phi(s, \lambda) \Psi(s, \lambda) + \Phi_\lambda(s, \lambda) \Psi + \Psi_\lambda(s, \lambda) \Phi(s, \lambda)] ds \notag\\
    & + k_1 \int_0^t [\bar{N}_{-\sigma} + \bar{N}_{-\sigma - \theta} + \Psi(s, \lambda)]\Phi(\lambda, s )ds\notag\\
    & \le \int_0^t[[ 2K_3 (\bar{B}_1 + \bar{B}_2 ) + k_1(\bar{N}_{-\sigma} + \bar{N}_{-\sigma - \theta} + \bar{B}_1)] \Phi(s, \lambda) + 2 k_3 \bar{B}_1 \Phi_\lambda(s, \lambda)] ds.
\end{align}
From the definition of $\Phi(t, \lambda,)$ in (\ref{defn_phi_psi}), we get
\begin{align}
    \Phi_\lambda (t, \lambda) = \int_0^\infty x \exp(\lambda x) |\phi(t, x)| dx.
\end{align}
From  (\ref{Uniq:eq_small_phi}), we obtain
\begin{align}\label{Uniq:phi_lambda}
\Phi_\lambda(t, \lambda) =& \int_0^t    \int_{x = 0}^\infty \left[x \exp(\lambda x) \sgn(\phi(s,x))\int_{z = 0}^\infty \int_{y=x}^\infty C(y, z) F(x,y|z) \right.\notag\\
&\left. [g_1(t, y) g_1(t, z) - g_2(t, y) g_2(t, z)] dy dz \right.\notag\\
    &\left. - \int_{y = 0}^\infty C(x, y) [g_1(t, x)g_1(t, y) - g_2(t, x) g_2(t,y)] dy\right] dx ds \notag\\
    = & \int_0^t \int_{z = 0}^\infty \int_{y = 0}^\infty \int_{x = 0}^y x \exp(\lambda x) C(y, z) F(x,y|z) G(s,y, z) dx dy dz ds \notag\\
    & - \int_0^t \int_{x=0}^\infty \int_{y = 0}^\infty x \exp(\lambda x) K(x, y) G(s, x, y) dy dx ds.
\end{align}
From first integral of (\ref{Uniq:phi_lambda}), we get
\begin{align}
    & \int_{z = 0}^\infty \int_{y = 0}^\infty \int_{x = 0}^y x \exp(\lambda x) C(y, z) F(x,y|z) G(s,y, z) dx dy dz ds \notag\\
    \le & k_1 k_2 \int_{z = 0}^\infty \int_{y = 0}^\infty \int_{x = 0}^y x \exp(\lambda x) (yz)^{- \sigma} y^{-\beta} G(s,y, z) dx dy dz ds \notag\\
    \le & k_1 k_2 \int_{z = 0}^\infty \int_{y = 0}^\infty  y \exp(\lambda y) (yz)^{- \sigma} y^{-\beta} G(s,y, z) dx dy dz ds \notag\\
    \le & k_1 k_2 \int_{z = 0}^\infty \int_{y = 0}^1   \exp(\lambda y) (y)^{1- \sigma- \beta} z^{-\sigma} G(s,y, z) dx dy dz ds \notag\\
    & + k_1 k_2 \int_{z = 0}^\infty \int_{y = 1}^\infty   \exp(\lambda y) (y)^{1- \sigma- \beta} z^{-\sigma} G(s,y, z) dx dy dz ds \notag\\
    \le & k_1 k_2 \int_{z = 0}^\infty \int_{y = 0}^1   \exp(\lambda y)  z^{-\sigma} [\phi(s, y)g_1(s, z) + g_2(s,y) \phi(s, z)] dx dy dz ds \notag\\
    & + k_1 k_2 \int_{z = 0}^\infty \int_{y = 1}^\infty   \exp(\lambda y) y z^{-\sigma}[\phi(s, y)g_1(s, z) + g_2(s,y) \phi(s, z)] dx dy dz ds \notag\\
    \le & k_1 k_2 [(\Phi(s, \lambda) \Psi(s, \lambda) + \Psi(s, \lambda) \Phi(s, \lambda)) + (\Phi_\lambda(s, \lambda) \Psi(s, \lambda) + \Psi_\lambda(s, \lambda) \Phi(s, \lambda))].
\end{align}
For the last integral of (\ref{Uniq:phi_lambda}), we get
\begin{align}
    & \int_{x=0}^\infty \int_{y = 0}^\infty x \exp(\lambda x) K(x, y) G(s, x, y) dy dx ds \notag\\
    =&k_1 \int_{y=0}^\infty \int_{x = 0}^1 x \exp(\lambda x) (xy)^{-\sigma} G(s, x, y) dy dx ds \notag\\
    & + k_1\int_{y=0}^\infty \int_{x = 1}^\infty x \exp(\lambda x) (xy)^{-\sigma} G(s, x, y) dy dx ds \notag\\
    = & k_1 \int_{y=0}^\infty \int_{x = 0}^1  \exp(\lambda x) y^{-\sigma} [\phi(s, x)g_1(s, y) + g_2(s,x) \phi(s, y)] dy dx ds \notag\\
    & + k_1\int_{y=0}^\infty \int_{x = 1}^\infty x \exp(\lambda x) y^{-\sigma} [\phi(s, x)g_1(s, y) + g_2(s,x) \phi(s, y)] dy dx ds \notag\\
    \le & k_1 [(\Phi(s, \lambda) \Psi(s, \lambda) + \Psi(s, \lambda) \Phi(s, \lambda)) + (\Phi_\lambda(s, \lambda) \Psi(s, \lambda) + \Psi_\lambda(s, \lambda) \Phi(s, \lambda))].
\end{align}
Therefore from (\ref{Uniq:phi_lambda}), we get
\begin{align}\label{ineq:2_phi_lambda}
    \Phi_\lambda(t, \lambda) & \le k_1 k_2 \int_0^t [2 \Phi(s, \lambda) \Psi(s, \lambda) + \Phi_\lambda(s, \lambda) \Psi(s, \lambda)+ + \Psi_\lambda(s, \lambda) \Phi(s, \lambda)]ds\notag\\
    & + \int_0^t k_1[2 \Phi(s, \lambda) \Psi(s, \lambda) + \Phi_\lambda(s, \lambda) \Psi(s, \lambda)+  \Psi_\lambda(s, \lambda) \Phi(s, \lambda)]ds \notag\\
    \le & 2 k_4\int_0^t [(2\bar{B}_1 + \bar{B}_2) \Phi(s, \lambda) + \bar{B}_1 \Phi_\lambda(s, \lambda)] ds, \text { where } k_4 = \max\{k_1 k_2, k_1\}.
\end{align}
Now, we recall a beautiful lemma to conclude the uniqueness proof.
\begin{lemma}\label{lema4.1} \textnormal{(See \cite{dubovskii1996existence}).}
Let the real-valued continuous function $w(t,\lambda)$ possesses continuous partial derivatives  $w_{\lambda}$ and $w_{\lambda \lambda}$ on $\mathcal{D}=\{ (t,\lambda): t\in[0,T], \lambda \in[0, \lambda_0] \}$. Further, assume that the real-valued functions  $\vartheta(t,\lambda)$, $\eta(\lambda)$, $\tau(t,\lambda)$ and $\rho(t,\lambda)$ and their  partial derivatives with respect to $\lambda$ are continuous on $\mathcal{D}$ and the function $w,w_{\lambda},\vartheta,\rho$ are nonnegative. Moreover, in $\mathcal{D}$  the following conditions hold
\[ w(t,\lambda) \leq \eta(\lambda) + \int_0^t (\rho(s,\lambda) v(s,\lambda) + \vartheta(s,\lambda)w_{\lambda}(s,\lambda) +  \tau(s,\lambda))ds\]
and
\[ w_{\lambda}(t,\lambda) \leq \eta_{\lambda}(\lambda) +\int_0^t \frac{\partial}{\partial \lambda}(\vartheta(s,\lambda)w_{\lambda}(s,\lambda) +\rho(s,\lambda)w(s,\lambda)+\tau(s,\lambda))ds. \]
Let $\mathcal{C}_0=\displaystyle\sup_{  \lambda \in[0, \lambda_0]} \eta, ~\mathcal{C}_1=\sup_\mathcal{D}\vartheta, ~ \mathcal{C}_2 = \sup_\mathcal{D} \rho$ and  $\mathcal{C}_3=\sup_\mathcal{D}\tau$. Then,
\[ w(t,\lambda)\leq \tfrac{\mathcal{C}_3}{\mathcal{C}_2} (\exp(\mathcal{C}_2t)-1) + \mathcal{C}_0 \exp (\mathcal{C}_2t)\]
in a region $\mathcal{R}\subset \mathcal{D}$:
\[\mathcal{ R}=\{(t,\lambda): \lambda \in[ \lambda_1-\mathcal{C}_1t,~ \lambda_0-\mathcal{C}_1t], ~ 0 \leq t \leq t^{'}<T^{'}, ~\lambda_1 \in (0,\lambda_0) \},\]
where $T^{'}= \min \{(\lambda_1/\mathcal{C}_1), T\}$.
\end{lemma}
Comparing (\ref{ineq_phi}) (\ref{ineq:2_phi_lambda}) with Lemma \ref{lema4.1} we get $\eta =0, \vartheta = 2 k_3 \bar{B}_1, \tau =0$ and $\rho = 2 k_4 (2 \bar{B}_1 + \bar{B}_2) \lambda + \bar{B}_1 + 2 k_3 (2 \bar{B}_1 + \bar{B}_2 + \bar{N}_{-\sigma} + \bar{N}_{-\sigma - \theta})$. Therefore, $\mathcal{C}_0 = 0 $ and $ \mathcal{C}_3 = 0$. Hence, $\Phi(t,\lambda) =0.$ By Lemma \ref{lema4.1} and applying  similar analysis in \cite{dubovskii1996existence}, we can conclude that $$\phi(t,x) =0, ~\text{i.e.},~ g_1(t, x)= g_2(t, x).$$ Hence, the uniqueness property is proved.
\end{proof}

\section{Concluding remarks}
\textcolor{black}{An extensive discussion of the existence and uniqueness result for the singular collision kernel model has been presented. To demonstrate that, first, we truncate the unbounded domain by compactly supported kernels. Next, we have shown the existence of a local solution by means of the Banach fixed point theorem. In the later stage, the existence of global solution in the complete domain is proved thanks to the Arzel$\grave{\text{a}}$-Ascoli theorem. Finally, with the help of Fubini's beautiful integral theorem, we also exhibit the conservation of mass property. Inclusion of singularity in the collisional kernel, the problem revisits earlier's existence results. We impose the least possible regulations over the initial data and fragmentation kernel to show the results. The singularity in the collision kernel includes several practical oriented kernels and could serve the development of many physical properties. It would be interesting if one can extend the existence result for $\sigma>\frac{1}{2}$.}\\ \\
\textbf{Acknowledgment:} Authors D. Ghosh and J. Kumar are thankful for funding support from Science and Engineering Research Board (SERB), Govt. of India (Grant number: EMR/2017/001514).

\bibliographystyle{plain}
\bibliography{existence_jp_dg}
\end{document}